\documentclass[onefignum,onetabnum]{siamart171218}
\usepackage{amsfonts}
\usepackage{graphicx}
\usepackage{epstopdf}
\usepackage{amsopn}
\usepackage{amstext,amssymb,amsmath,amsbsy}
\usepackage{hyperref}
\usepackage{amscd}
\usepackage{amsfonts}
\usepackage{indentfirst}
\usepackage{verbatim}
\usepackage{amsmath}
\usepackage{enumerate}
\usepackage{graphicx}
\usepackage{color}
\usepackage[OT1]{fontenc}
\usepackage[latin1]{inputenc}
\usepackage[english]{babel}
\usepackage{amssymb}
\usepackage{subfig}
\usepackage{algorithm}
\usepackage{algpseudocode}

\ifpdf
  \DeclareGraphicsExtensions{.eps,.pdf,.png,.jpg}
\else
  \DeclareGraphicsExtensions{.eps}
\fi

\newsiamremark{remark}{Remark}
\newsiamremark{hypothesis}{Hypothesis}
\crefname{hypothesis}{Hypothesis}{Hypotheses}
\newsiamthm{claim}{Claim}
\headers{Linearization of the kinematic inverse problem}{M. V. Klibanov, T. T. Le and L. H. Nguyen }
\newcommand{\R}{\mathbb{R}}

\newcommand{\x}{{\bf x}}

\newtheorem{Proposition}{Proposition}[section]

\newtheorem{Remark}{Remark}[section]

\newtheorem{Assumption}{Assumption}[section]
\newtheorem*{assumption*}{Assumption}

\newtheorem{problem}{Problem}[section]

\makeatletter
\newcommand*{\addFileDependency}
[1]{
  \typeout{(#1)}
  \@addtofilelist{#1}
  \IfFileExists{#1}{}{\typeout{No file #1.}}
}
\makeatother
\newcommand*{\myexternaldocument}
[1]{
    \externaldocument{#1}
    \addFileDependency{#1.tex}
    \addFileDependency{#1.aux}
}
\myexternaldocument{ex_supplement}
\input{tcilatex}
\begin{document}

\title{Convergent numerical method for a linearized travel time tomography
problem with incomplete data \thanks{
{\textbf{Funding.} The work was supported by US Army Research Laboratory and
US Army Research Office grant W911NF-19-1-0044.}} }
\author{Michael V. Klibanov\thanks{%
Department of Mathematics and Statistics, University of North Carolina at
Charlotte, Charlotte, NC 28223, USA, mklibanv@uncc.edu (corresponding
author), tle55@uncc.edu, loc.nguyen@uncc.edu} \and Thuy T. Le\footnotemark[2]
\and Loc H. Nguyen\footnotemark[2] }
\maketitle

\begin{abstract}
We propose a new numerical method to solve the linearized problem of travel
time tomography with incomplete data. Our method is based on the technique
of the truncation of the Fourier series with respect to a special basis of $%
L^{2}$. This way we derive a boundary value problem for a system of coupled
partial differential equations (PDEs) of the first order. This system is
solved by the quasi-reversibility method. Hence, the spatially dependent
Fourier coefficients of the solution to the linearized Eikonal equation are
obtained. The convergence of this method is established. Numerical results
for highly noisy data are presented.
\end{abstract}

% REQUIRED

% REQUIRED

% REQUIRED

\section{Introduction}

\label{sec 1}

In this paper we develop a new numerical method for the linearized Travel
Time Tomography Problem (TTTP) for the $d-$D case. Our data are both
non-redundant and incomplete. Using a discrete Carleman estimate, we
establish the convergence of our method. In addition, we provide results of
numerical experiments in the 2D case. In particular, we demonstrate that our
method provides good accuracy of images of complicated objects with 5\%
noise in the data. Furthermore, a satisfactory accuracy of images is
demonstrated even for very high levels of noise between 30\% and 120\%.

In fact, both the idea of our method and sources/detectors configuration are
close to those of our recent works \cite%
{KlibanovNguyen:ip2019,KlibanovAlexeyNguyen:SISC2019}. However, our case is
substantially more difficult one since the waves in our case propagate along
geodesic lines, rather than a radiation propagating along straight lines in 
\cite{KlibanovNguyen:ip2019,KlibanovAlexeyNguyen:SISC2019}. Still, although
we formulate here results related to the convergence of our method, we do
not prove them. The reason is that, as it turns out, proofs are very similar
to those in \cite{KlibanovAlexeyNguyen:SISC2019}. In other words,
surprisingly, the analytical apparatus of the convergence theory developed
in \cite{KlibanovAlexeyNguyen:SISC2019} works well for the problem
considered in this paper.

In the isotropic case of acoustic/seismic waves propagation, the travel time
tomography problem (TTTP) is the problem of the recovery of the spatially
distributed speed of propagation of acoustic/seismic waves from the first
times of arrival of those waves. In the electromagnetic case this is the
problem of the recovery of the spatially distributed dielectric constant
from those times. Another name for the TTTP is inverse kinematic problem
(IKP). Waves are originated by some sources located either at the boundary
of the closed bounded domain of interest or outside of this domain. Times of
first arrival from those sources are measured on a part of the boundary of
that domain. The TTTP has well known applications in Geophysics, see, e.g.
the book of Romanov \cite[Chapter 3]{Romanov:VNU1986}.

The history of the TTTP has started 114 years ago. The pioneering papers
about the solution of the 1D TTTP were published by Herglotz \cite%
{Herglotz:zmp1905} (1905) and then by Wiechert and Zoeppritz \cite%
{WiechertZoeppritz:gwg1097} (1907). Their method is described in the book of
Romanov \cite[Section 3 of Chapter 3]{Romanov:VNU1986}. It was discovered
recently that, in addition to Geophysics, the IKP has applications in the
phaseless inverse scattering problem \cite{KlibanovLocKejia:apnum2016,
KlibanovRomanov:SIAMam2016, Romanov:dm2017}.

The next natural question after the classical 1D case of \cite%
{Herglotz:zmp1905, WiechertZoeppritz:gwg1097} was about $2$ and $3$
dimensional cases. The first uniqueness and Lipschitz stability result for
the 2D case was obtained by Mukhometov \cite{Mukhometov:smd1977}, also see 
\cite{BernsteinGerver:dan1978, PestovUhlmann:am2005}. Next, these results
were obtained by Mukhometov and Romanov for the 3D case in \cite%
{MukhometovRomanov:das1978, Romanov:VNU1986}. We also refer to the work of
Stefanov, Uhlmann and Vasy \cite{StefanovUhlmannVasy:arxiv2017} for a more
recent publication for the 3D case. As to the numerical methods for the
inverse kinematic problem, we refer to \cite{SchrHoderSchuster:ip2016} for
the 2D case and to \cite{ZhaoZhong:SIIS2019} for the 3D case.

In all past publication about the IKP, the data are redundant in the 3D case
and complete in both 2D and 3D cases. In two recent works of the first
author \cite{Klibanov:jiip2019, Klibanov:IPI2019} two globally convergent
numerical methods for the 3D TTTP with non redundant incomplete data were
developed.

Along with the full IKP, a significant applied interest is also in a
linearized IKP, see \cite[Chapter 3]{Romanov:VNU1986}. Let $c$ be the speed
of sound. Denote $\mathbf{n}=1/c$ the refractive index. To linearize, one
should assume that $\mathbf{n}=\mathbf{n}_{0}+\mathbf{n}_{1},$ where $%
\mathbf{n}_{0}$ is the known background function and $\mathbf{n}_{1}$ with $%
\left\vert \mathbf{n}_{1}\right\vert \ll \mathbf{n}_{0}$ is its unknown
perturbation, which is the subject to the solution of the linearized TTTP.
Thus, one assumes that the refractive index is basically known, whereas its
small perturbation $\mathbf{n}_{1}$ is unknown. This problem is also called
the \emph{geodesic X-ray transform problem}. The Lipschitz stability and
uniqueness theorem for this problem in the isotropic case was first obtained
in \cite{Romanov:dan1978}, see Theorem 3.2 in Section 4 of Chapter 3 of \cite%
{Romanov:VNU1986}. In the non isotropic case this problem was studied in 
\cite{StefanovUhlmannVasy:jam2018}. In \cite{Monard:SIAM2014} numerical
studies of this problem in the isotropic case were performed.

In our derivation, we end up with an over determined boundary value problem
for a system of coupled linear PDEs of the first order. It is well known
that the quasi-reversibility method is an effective tool for numerical
solutions of over determined boundary value problems for PDEs. Lattès and
Lions \cite{LattesLions:e1969} were the first ones who have proposed the
quasi-reversibility method. This technique was developed further in, e.g. 
\cite{BourgeoisDarde:ip2010, BourgeoisPonomarevDarde:ipi2019,
Isakov:bookSpringer2017,Klibanov:anm2015,KlibanovNguyen:ip2019,
LiNguyen:IPSE2019,KlibanovAlexeyNguyen:SISC2019}. In particular, it was
shown in \cite{Klibanov:anm2015} that while it is rather easy to prove,
using Riesz theorem, the existence and uniqueness of the minimizer of a
certain functional related to this method, the proof of convergence of those
minimizers to the correct solution requires a stronger tool of Carleman
estimates.

Another important feature of this paper is a special orthonormal basis in
the space $L^{2}\left( -\overline{\alpha },\overline{\alpha }\right) ,$
where $\overline{\alpha }>0$ is a certain number. The functions of this
basis depend only on the position of the point source. This basis was first
introduced in \cite{Klibanov:jiip2017} and was further used in \cite%
{Klibanov:jiip2019,
Klibanov:IPI2019,KlibanovNguyen:ip2019,KlibanovAlexeyNguyen:SISC2019}. Just
like in our previous publications \cite{Klibanov:jiip2019,
Klibanov:IPI2019,KlibanovNguyen:ip2019,KlibanovAlexeyNguyen:SISC2019}, we
use here an \emph{approximate mathematical model}. More precisely, we assume
that a certain function associated with the solution of the governing
linearized Eikonal equation can be represented via a truncated Fourier
series with respect to this basis. This assumption forms the first element
of that model. The second element is that we assume that the first
derivatives with respect to all variables, except of one, are written via
finite differences and the step size of these finite differences is bounded
from the below by a positive number $h_{0}>0$.

We do not prove convergence as the number $N$ of terms in that truncated
series tends to infinity and the lower bound for the grid step $h_{0}$ size
tends to zero. Thus, we come up with a semi-finite dimensional approximate
mathematical model. We point out that similar approximate mathematical
models are used quite often in studies of numerical methods for inverse
problems by other authors, and numerical results are usually encouraging,
see, e.g. \cite%
{GuillementNovikov:ipse2013,Kabanikhin:n1988,KabanikhinSatybaevShishlenin:svp2005,KabanikhinSabelfeldNovikovShishlenin:jiip2015}%
. Just as ourselves, proofs of convergence results in such cases when, e.g. $%
N\rightarrow \infty ,h_{0}\rightarrow 0$ are usually not conducted since
they are very challenging tasks due to the ill-posed nature of inverse
problems.

The paper is organized as follows. In Section \ref{sec statement}, we
formulate the inverse problem. Next, in Section \ref{sec truncation}, we
introduce the truncation technique and our numerical method. Then, in
Section \ref{sec quasi}, we recall the quasi-reversibility method and its
convergence in the case of partial finite differences. 
In Section \ref{Sec num} we present the implementation and numerical
results. Finally, Section \ref{sec conclusion} is for concluding remarks.

\section{The linearization}

\label{sec statement} Let $d\geq 2$ be the spatial dimension. Let $R>1$ and $%
0<a<b$. Set 
\begin{equation}
\Omega =(-R,R)^{d-1}\times (a,b)\subset \mathbb{R}^{d}.  \label{1}
\end{equation}%
Let $\mathbf{c}_{0}(\mathbf{x})=\mathbf{n}_{0}^{2}(\mathbf{x})$, $\mathbf{x}%
\in \Omega $ where $\mathbf{n}_{0}$ is the refractive index of the
background. Assume that $\mathbf{c}_{0}=\mathbf{n}_{0}^{2}=1$ on $\mathbb{R}%
^{d}\setminus \Omega .$ For any two points $\mathbf{x}_{1}$ and $\mathbf{x}%
_{2}$ in $\mathbb{R}^{d}$, define the geodesic line generated by $\mathbf{n}%
_{0}$ connecting $\mathbf{x}_{1}$ and $\mathbf{x}_{2}$ as: 
\begin{multline}
\Gamma (\mathbf{x}_{1},\mathbf{x}_{2})=\mathrm{argmin}\Big\{\int_{\gamma }%
\mathbf{n}_{0}(\mathbf{\xi })d\sigma (\mathbf{\xi })\mbox{ where }\gamma
:[0,1]\rightarrow \mathbb{R}^{d} \\
\mbox{ is a smooth map with }\gamma (0)=\mathbf{x}_{1},\gamma (1)=\mathbf{x}%
_{2}\Big\}.  \label{2}
\end{multline}%
Given the refractive index $\mathbf{n}_{0}$, the geodesic line $\Gamma (%
\mathbf{x}_{1},\mathbf{x}_{2})$ is the curve connecting points $\mathbf{x}%
_{1}$ and $\mathbf{x}_{2}$ and such that the travel time along $\Gamma (%
\mathbf{x}_{1},\mathbf{x}_{2})$ is minimal. The travel time is the integral
in \eqref{2}. If $\mathbf{n}_{0}\equiv 1$, then $\Gamma (\mathbf{x}_{1},%
\mathbf{x}_{2})$ is the line segment connecting these two points.

Introduce the line of sources $L_{\mathrm{sc}}$ located on the $x_{1}$-axis
as 
\begin{equation}
L_{\mathrm{sc}}=[-\overline{\alpha },\overline{\alpha }]\times \{(0,0,\dots
,0)\},  \label{def Lsc}
\end{equation}%
where $\overline{\alpha }$ is a fixed positive number. For each source
position $\mathbf{x}_{\alpha } = (\alpha, 0, \dots, 0) \in L_{\mathrm{sc}}$,
the function 
\begin{equation}
u_{0}(\mathbf{x},\mathbf{x}_{\alpha })=\int_{\Gamma (\mathbf{x,x}_{\alpha })}%
\mathbf{n}_{0}(\mathbf{\xi })d\sigma (\mathbf{\xi })\quad \mathbf{x}\in 
\mathbb{R}^{d}  \label{def u}
\end{equation}%
is the travel time of the wave from $\mathbf{x}_{\alpha }$ to $\mathbf{x}$.

\begin{Assumption}[regularity of geodesic lines]
We assume everywhere in this paper that the geodesic lines are regular in
the following sense: for each point $\mathbf{x}$ of the closed domain $%
\overline{\Omega }$ and for each point $\mathbf{x}_{\alpha }$ of the line of
sources $L_{\mathrm{sc}}$ there exists a single geodesic line $\Gamma (%
\mathbf{x,x}_{\alpha })$ connecting them.
\end{Assumption}

For each $\alpha \in (-\overline \alpha, \overline \alpha),$ define 
\begin{align*}
\partial \Omega _{\alpha }^{-} &= \{\mathbf{x}\in \partial \Omega :\nabla
u_{0}(\mathbf{x},\mathbf{x}_{\alpha })\cdot n(\mathbf{x})\leq 0\}, \\
\partial \Omega _{\alpha }^{+} &= \{\mathbf{x}\in \partial \Omega :\nabla
u_{0}(\mathbf{x},\mathbf{x}_{\alpha })\cdot n(\mathbf{x})>0\}
\end{align*}%
where $\mathbf{x}_\alpha = (\alpha, 0, \dots, 0).$ Let $p:\mathbb{R}%
^{d}\rightarrow \mathbb{R}$ be a function compactly supported in $\Omega $.
For each $\mathbf{x}_{\alpha }\in L_{\mathrm{sc}}$, let $u(\mathbf{x},%
\mathbf{x}_{\alpha })$ be the solution to 
\begin{equation}
\left\{ 
\begin{array}{rcll}
\nabla u_{0}(\mathbf{x},\mathbf{x}_{\alpha })\cdot \nabla u(\mathbf{x},%
\mathbf{x}_{\alpha }) & = & p(\mathbf{x}), & \mathbf{x}\in \Omega , \\ 
u(\mathbf{x},\mathbf{x}_{\alpha }) & = & 0, & \mathbf{x}\in \Omega ^{-}.%
\end{array}%
\right.  \label{main eqn}
\end{equation}%
The aim of this paper is to solve the following inverse problem:

\begin{problem}[linearized travel time tomography problem]
Given the data 
\begin{equation}
f(\mathbf{x},\mathbf{x}_{\alpha })=u(\mathbf{x},\mathbf{x}_{\alpha }),\quad 
\mathbf{x}\in \partial \Omega ^{+},\mathbf{x}_{\alpha }\in L_{\mathrm{sc}},
\label{5}
\end{equation}
determine the function $p(\mathbf{x}),$ $\mathbf{x}\in \Omega .$ \label%
{problem isp}
\end{problem}

\begin{remark}
The data $f(\mathbf{x},\mathbf{x}_{\alpha })$ are non-redundant ones.
Indeed, the source $\mathbf{x}_{\alpha }$ $\in L_{\mathrm{sc}}$ depends on
one variable and $\mathbf{x}\in \partial \Omega ^{+}$ depends on $d-1$
variables. Hence the function $f(\mathbf{x},\mathbf{x}_{\alpha })$ depends
on $d$ variables, so does the target function $p(\mathbf{x}).$ \label{rem
nonredundant}
\end{remark}

Problem \ref{problem isp} arises from the highly nonlinear and severely
ill-posed inverse kinematic problem. Assume that $\mathbf{c}(\mathbf{x})=%
\mathbf{n}^{2}(\mathbf{x})$ contains a perturbation term of the background
function $\mathbf{c}_{0}(\mathbf{x})=\mathbf{n}_{0}^{2}(\mathbf{x})$. In
other words, 
\begin{equation}
\mathbf{c}(\mathbf{x})=\mathbf{c}_{0}(\mathbf{x})+2\epsilon p(\mathbf{x}%
)\quad \mathbf{x}\in \mathbb{R}^{d}  \label{eqn pur}
\end{equation}%
for a small number $\epsilon >0.$ Denote by 
\begin{equation*}
u_{\mathbf{n}}(\mathbf{x},\mathbf{x}_{\alpha })=\int_{\Gamma _{\mathbf{n}}(%
\mathbf{x},\mathbf{x}_{\alpha })}\mathbf{n}(\xi )d\sigma (\xi )
\end{equation*}%
the travel time from $\mathbf{x}_{\alpha }\in L_{\mathrm{sc}}$ to $\mathbf{x}%
\in \Omega $, where $\Gamma _{\mathbf{n}}(\mathbf{x},\mathbf{x}_{\alpha })$
is the geodesic line generated by the function $\mathbf{n}$. Then, it is
well-known \cite{Romanov:VNU1986} that $u(\mathbf{x},\mathbf{x}_{\alpha })$
satisfies the Eikonal equation 
\begin{equation}
|\nabla u_{\mathbf{n}}(\mathbf{x},\mathbf{x}_{\alpha })|^{2}=\mathbf{c}(%
\mathbf{x})\quad \mathbf{x}\in \Omega ,\mathbf{x}_{\alpha }\in L_{\mathrm{sc}%
}.  \label{eqn Eik}
\end{equation}%
The inverse kinematic problem is to determine the function $\mathbf{c}$ from
the measurement of $u(\mathbf{x},\mathbf{x}_{\alpha })$ for all $\mathbf{x}%
\in \partial \Omega ^{+}$ and $\mathbf{x}_{\alpha }\in L_{\mathrm{sc}}$. Let 
$u_{0}(\mathbf{x},\mathbf{x}_{\alpha })$ be the travel time function
corresponding to the background $\mathbf{c}_{0}$. Then, one has 
\begin{equation}
|\nabla u_{0}(\mathbf{x},\mathbf{x}_{\alpha })|^{2}=\mathbf{c}_{0}(\mathbf{x}%
)\quad \mathbf{x}\in \Omega ,\mathbf{x}_{\alpha }\in L_{\mathrm{sc}}.
\label{eqn Eik0}
\end{equation}%
Due to \eqref{eqn pur} we represent $\nabla u_{\mathbf{n}}(\mathbf{x},%
\mathbf{x}_{\alpha })$ as $\nabla u_{\mathbf{n}}(\mathbf{x},\mathbf{x}%
_{\alpha })=\nabla u_{0}(\mathbf{x},\mathbf{x}_{\alpha })+\epsilon \nabla
u^{\left( 1\right) }(\mathbf{x},\mathbf{x}_{\alpha }).$ Hence, ignoring the
term with $\epsilon ^{2},$ we obtain 
\begin{equation*}
|\nabla u_{\mathbf{n}}(\mathbf{x},\mathbf{x}_{\alpha })|^{2}\approx |\nabla
u_{0}(\mathbf{x},\mathbf{x}_{\alpha })|^{2}+2\epsilon \nabla u_{0}(\mathbf{x}%
,\mathbf{x}_{\alpha })\nabla u^{\left( 1\right) }(\mathbf{x},\mathbf{x}%
_{\alpha }).
\end{equation*}%
Denoting $u^{\left( 1\right) }:=u,$ we obtain 
\begin{equation*}
\nabla u_{0}(\mathbf{x},\mathbf{x}_{\alpha })\cdot \nabla u(\mathbf{x},%
\mathbf{x}_{\alpha })=p(\mathbf{x}).
\end{equation*}%
Thus, the inverse source problem under consideration is the
\textquotedblleft linearization" of the nonlinear kinematic inverse problem.

%On the other hand, Problem \ref{problem isp} can be consider as the problem of inverting the geodesic Radon transform.
%Assume that $q$

Note that since the function $p$ is compactly supported in $\Omega ,$ then $%
\mathbf{c}=\mathbf{c}_{0}=1$ in $\mathbb{R}^{d}\setminus \Omega $. This
implies that $u_{\mathbf{n}}(\mathbf{x},\mathbf{x}_{\alpha })=u_{0}(\mathbf{x%
},\mathbf{x}_{\alpha })=|\mathbf{x}-\mathbf{x}_{\alpha }|,$ for all $\mathbf{%
x}\in \partial \Omega _{\alpha }^{-}.$ Hence, we set %$\begin{equation*}
$u(\mathbf{x},\mathbf{x}_{\alpha })=0$ for all $\mathbf{x}\in \partial
\Omega _{\alpha }^{-}. $ %\end{equation*}

From now on, to separate the coordinate number $d$ of the point $\mathbf{x}$%
, we write $\mathbf{x}=(x_{1},\dots ,x_{d-1},z)$. The transport equation in %
\eqref{main eqn} is read as 
\begin{equation}
\partial _{z}u_{0}(\mathbf{x},\mathbf{x}_{\alpha })\partial _{z}u(\mathbf{x},%
\mathbf{x}_{\alpha })+\sum_{i=1}^{d-1}\partial _{x_{i}}u_{0}(\mathbf{x},%
\mathbf{x}_{\alpha })\partial _{x_{i}}u(\mathbf{x},\mathbf{x}_{\alpha })=p(%
\mathbf{x})  \label{2.10}
\end{equation}%
for all $\mathbf{x}\in \Omega ,$ $\mathbf{x}_{\alpha }\in L_{\mathrm{sc}}$.

\section{A boundary value problem for a system of coupled PDEs of the first
order}

\label{sec truncation}

This section aims to derive a system of partial differential equations,
which can be stably solved by the quasi-reversibility method in the
semi-finite difference scheme. The solution of this system yields the
desired numerical solution to Problem \ref{problem isp}.

%\subsection{Preliminaries}

We will employ a special basis of $L^2(-\overline \alpha, \overline \alpha)$
where $2\overline \alpha$ is the length of the line of source $L_{\mathrm{sc}%
},$ see \eqref{def Lsc}. For each $n = 1, 2, \cdots,$ let $\phi_n(\alpha) =
\alpha^{n - 1}\exp(\alpha)$. The set $\{\phi_n\}_{n = 1}^{\infty}$ is
complete in $L^2(-\overline \alpha, \overline \alpha).$ Applying the
Gram-Schmidt orthonormalization process to this set, we obtain a basis of $%
L^2(-\overline \alpha, \overline \alpha)$, named as $\{\Psi_n\}_{n =
1}^{\infty}$. We have the proposition

\begin{Proposition}[see \protect\cite{Klibanov:jiip2017}]
The basis $\{\Psi_n\}_{n = 1}^{\infty}$ satisfies the following properties:

\begin{enumerate}
\item $\Psi_n$ is not identically zero for all $n \geq 1$, %

\item For all $m,n\geq 1$ 
\begin{equation*}
s_{mn}=\int_{-\overline{\alpha }}^{\overline{\alpha }}\Psi _{n}^{\prime
}(\alpha )\Psi _{m}(\alpha )d\alpha =\left\{ 
\begin{array}{ll}
1 & \mbox{if }m=n, \\ 
0 & \mbox{if }n<m.%
\end{array}%
\right.
\end{equation*}%
As a result, for all integer $N>1$, the matrix $S_{N}=(s_{mn})_{m,n=1}^{N}$,
is invertible.
\end{enumerate}

\label{prop MK}
\end{Proposition}

\begin{Remark}
The basis $\{\Psi _{n}\}_{n=1}^{\infty }$ was first introduced in \cite%
{Klibanov:jiip2017}. Then, this basis was successfully used to solve several
important inverse problems, including the inverse source problem for
Helmholtz equations \cite{NguyenLiKlibanov:IPI2019}, inverse X-ray
tomographic problem in incomplete data \cite{KlibanovNguyen:ip2019} and the
nonlinear inverse problem of electrical impedance tomography with the
so-called restricted Dirichlet-to-Neumann map data, see \cite%
{KlibanovLiZhang:ip2019}, the inverse problem of computing the initial
condition of nonlinear parabolic equations \cite{LeNguyen:arxiv2019}.
\end{Remark}

%\subsection{A system of transport equations}

We now derive an important system for Fourier coefficients of the function 
\begin{equation}
w(\mathbf{x},\mathbf{x}_{\alpha })=u(\mathbf{x},\mathbf{x}_{\alpha
})\partial _{z}u_{0}(\mathbf{x},\mathbf{x}_{\alpha })\quad \mathbf{x}\in
\Omega ,\mathbf{x}_{\alpha }\in L_{\mathrm{sc}}  \label{3.1}
\end{equation}%
with respect to the basis in Proposition \ref{prop MK}. Differentiate %
\eqref{2.10} with respect to $\alpha $. We obtain 
\begin{equation}
\frac{\partial }{\partial \alpha }\Big[\partial _{z}u_{0}(\mathbf{x},\mathbf{%
x}_{\alpha })\partial _{z}u(\mathbf{x},\mathbf{x}_{\alpha
})+\sum_{i=1}^{d-1}\partial _{x_{i}}u_{0}(\mathbf{x},\mathbf{x}_{\alpha
})\partial _{x_{i}}u(\mathbf{x},\mathbf{x}_{\alpha })\Big]=0
\label{v and valpha}
\end{equation}%
for all $\mathbf{x}\in \Omega ,$ $\mathbf{x}_{\alpha }\in L_{\mathrm{sc}}.$
From now on, we impose the following condition.

\begin{Assumption}[Monotonicity condition in the $z$-direction]
The traveling time function $u_{0}$, defined in \eqref{def u} with $\mathbf{n%
}$ replaced by $\mathbf{n}_{0}$, is strictly increasing with respect to $z$.
In other words, 
\begin{equation*}
\partial _{z}u_{0}(\mathbf{x},\mathbf{x}_{\alpha })=\frac{\partial u_{0}(%
\mathbf{x},\mathbf{x}_{\alpha })}{\partial z}>0
\end{equation*}%
for all $\mathbf{x}=(x_{1},\dots ,x_{d-1},z)\in \Omega $ and for all $%
\mathbf{x}_{\alpha }\in L_{\mathrm{sc}}.$ \label{cond uy}
\end{Assumption}

Assumption \ref{cond uy} means that the higher in the $z$-direction, the
longer the traveling time is. A sufficient condition for Assumption \ref%
{cond uy} to be true is formulated in \eqref{4} of Lemma \ref{lemma 3.1}. A similar
monotonicity condition can be found in formulas (3.24) and (3.24$^{\prime })$
of section 2 of chapter 3 of the book \cite{Romanov:VNU1986}. Also, a
similar condition was imposed in originating works \ for the 1D problem of
Herglotz and Wiechert and Zoeppritz \cite%
{Herglotz:zmp1905,WiechertZoeppritz:gwg1097}: see section 3 of chapter 3 of 
\cite{Romanov:VNU1986}. Besides, figures 5 and 10 of \cite{Volgyesi:PPCE1982}
justify this condition from the geophysical standpoint. Although Lemma 3.1
is proven in \cite{Klibanov:IPI2019} only in the 3D case, the proof in the $%
d-$D case is very similar and, therefore, avoided.

\begin{lemma}[\protect\cite{Klibanov:IPI2019}]
%\textbf{Lemma 3.1 }\cite{Klibanov:IPI2019}.
Let conditions (\ref{1}) and (\ref{def Lsc}) hold. Also, assume that $%
\mathbf{c}_{0}\in C^{2}(\mathbb{R}^{d}),$ $\mathbf{c}_{0}\geq m_{0}$ for
some positive constant $m_{0}$, $\mathbf{c}_{0}(\mathbf{x})=1$ for all $z<a$
\begin{equation}
\partial _{z}c_{0}\left( \mathbf{x}\right) \geq 0\quad 
\mbox{for all }%
\mathbf{x}\in \overline{\Omega }.  \label{4}
\end{equation%
}%
Then, 
\begin{equation*}
\partial _{z}u_{0}(\mathbf{x},\mathbf{x}%
_{\alpha })\geq \frac{a}{\sqrt{%
a^{2}+2}}\quad \mbox{for all }\mathbf{x}%
\in \overline{\Omega },\alpha \in %
\left[ -\overline{\alpha },\overline{%
\alpha }\right] .
\end{equation*}
\label{lemma 3.1}
\end{lemma}

Consider a new function $w(\mathbf{x},\mathbf{x}_{\alpha }),$ 
\begin{equation}
w(\mathbf{x},\mathbf{x}_{\alpha })=u(\mathbf{x},\alpha )\partial _{z}u_{0}(%
\mathbf{x},\mathbf{x}_{\alpha })\quad \mathbf{x}\in \Omega ,\mathbf{x}%
_{\alpha }\in L_{\mathrm{sc}}.  \label{change variable}
\end{equation}%
We have 
\begin{align}
\partial _{z}u_{0}(\mathbf{x},\mathbf{x}_{\alpha })\partial _{z}u(\mathbf{x},%
\mathbf{x}_{\alpha })& =\partial _{z}w(\mathbf{x},\mathbf{x}_{\alpha })-u(%
\mathbf{x},\alpha )\partial _{zz}u_{0}(\mathbf{x},\mathbf{x}_{\alpha })
\label{3.3} \\
& =\partial _{z}w(\mathbf{x},\mathbf{x}_{\alpha })-w(\mathbf{x},\alpha )%
\frac{\partial _{zz}u_{0}(\mathbf{x},\mathbf{x}_{\alpha })}{\partial
_{z}u_{0}(\mathbf{x},\mathbf{x}_{\alpha })}.  \notag
\end{align}%
and for $i=1,\dots ,d-1$ 
\begin{align}
\partial _{x_{i}}u(\mathbf{x},\mathbf{x}_{\alpha })& =\frac{\partial }{%
\partial x_{i}}\left( \frac{w(\mathbf{x},\mathbf{x}_{\alpha })}{\partial
_{z}u_{0}(\mathbf{x},\mathbf{x}_{\alpha })}\right)  \label{3.4444} \\
& =\frac{\partial _{x_{i}}w(\mathbf{x},\mathbf{x}_{\alpha })\partial
_{z}u_{0}(\mathbf{x},\mathbf{x}_{\alpha })-w(\mathbf{x},\mathbf{x}_{\alpha
})\partial _{zx_{i}}u_{0}(\mathbf{x},\mathbf{x}_{\alpha })}{(\partial
_{z}u_{0}(\mathbf{x},\mathbf{x}_{\alpha }))^{2}}  \notag
\end{align}%
for all $\mathbf{x}\in \Omega ,\mathbf{x}_{\alpha }\in L_{\mathrm{sc}}.$
Combining \eqref{v and valpha}, \eqref{3.3} and \eqref{3.4444}, we obtain 
\begin{multline}
\frac{\partial }{\partial \alpha }\Big[\partial _{z}w(\mathbf{x},\mathbf{x}%
_{\alpha })-w(\mathbf{x},\mathbf{x}_{\alpha })\frac{\partial _{zz}u_{0}(%
\mathbf{x},\mathbf{x}_{\alpha })}{\partial _{z}u_{0}(\mathbf{x},\mathbf{x}%
_{\alpha })}  \label{56} \\
+\sum_{i=1}^{d-1}\frac{\partial _{x_{i}}w(\mathbf{x},\mathbf{x}_{\alpha
})\partial _{z}u_{0}(\mathbf{x},\mathbf{x}_{\alpha })-w(\mathbf{x},\mathbf{x}%
_{\alpha })\partial _{zx_{i}}u_{0}(\mathbf{x},\mathbf{x}_{\alpha })}{%
(\partial _{z}u_{0}(\mathbf{x},\mathbf{x}_{\alpha }))^{2}}\partial
_{x_{i}}u_{0}(\mathbf{x},\mathbf{x}_{\alpha })\Big]=0.
\end{multline}%
This is equivalent to 
\begin{multline}
\partial _{\alpha z}w(\mathbf{x},\mathbf{x}_{\alpha })-\frac{\partial
_{zz}u_{0}(\mathbf{x},\mathbf{x}_{\alpha })}{\partial _{z}u_{0}(\mathbf{x},%
\mathbf{x}_{\alpha })}\partial _{\alpha }w(\mathbf{x},\mathbf{x}_{\alpha })-%
\frac{\partial }{\partial \alpha }\left( \frac{\partial _{zz}u_{0}(\mathbf{x}%
,\mathbf{x}_{\alpha })}{\partial _{z}u_{0}(\mathbf{x},\mathbf{x}_{\alpha })}%
\right) w(\mathbf{x},\mathbf{x}_{\alpha }) \\
+\sum_{i=1}^{d-1}\Big[\frac{\partial _{x_{i}}u_{0}(\mathbf{x},\mathbf{x}%
_{\alpha })}{\partial _{z}u_{0}(\mathbf{x},\mathbf{x}_{\alpha })}\partial
_{\alpha x_{i}}w(\mathbf{x},\mathbf{x}_{\alpha })+\frac{\partial }{\partial
\alpha }\left( \frac{\partial _{x_{i}}u_{0}(\mathbf{x},\mathbf{x}_{\alpha })%
}{\partial _{z}u_{0}(\mathbf{x},\mathbf{x}_{\alpha })}\right) \partial
_{x_{i}}w(\mathbf{x},\mathbf{x}_{\alpha }) \\
-\frac{\partial _{zx_{i}}u_{0}(\mathbf{x},\mathbf{x}_{\alpha })\partial
_{x_{i}}u_{0}(\mathbf{x},\mathbf{x}_{\alpha })}{(\partial _{z}u_{0}(\mathbf{x%
},\mathbf{x}_{\alpha }))^{2}}\partial _{\alpha }w(\mathbf{x},\mathbf{x}%
_{\alpha }) \\
-\frac{\partial }{\partial \alpha }\left( \frac{\partial _{zx_{i}}u_{0}(%
\mathbf{x},\mathbf{x}_{\alpha })\partial _{x_{i}}u_{0}(\mathbf{x},\mathbf{x}%
_{\alpha })}{(\partial _{z}u_{0}(\mathbf{x},\mathbf{x}_{\alpha }))^{2}}%
\right) w(\mathbf{x},\mathbf{x}_{\alpha })\Big]=0.  \label{3.4}
\end{multline}

We recall now the orthonormal basis $\{\Psi _{n}\}_{n=1}^{\infty }$
constructed at the beginning of this section. For each $\mathbf{x}\in \Omega 
$ and for all $\mathbf{x}_{\alpha }\in L_{\mathrm{sc}}$, we write 
\begin{equation}
w(\mathbf{x},\mathbf{x}_{\alpha })=\sum_{n=1}^{\infty }w_{n}(\mathbf{x})\Psi
_{n}(\alpha )\approx \sum_{n=1}^{N}w_{n}(\mathbf{x})\Psi _{n}(\alpha ),
\label{Fourier w}
\end{equation}%
\begin{equation}
w_{n}(\mathbf{x})=\int_{-\overline{\alpha }}^{\overline{\alpha }}w(\mathbf{x}%
,\mathbf{x}_{\alpha })\Psi _{n}(\alpha )d\alpha .
\label{Fourier Coefficient}
\end{equation}%
The \textquotedblleft cut-off" number $N$ is chosen numerically. We discuss
the choice of $N$ in more details in Section \ref{Sec num}. Following our
approximate mathematical model introduced in Section 1, \ we assume that the
approximation $\approx $ in \eqref{Fourier w} is an equality as well as 
\begin{equation}
\partial _{\alpha }w(\mathbf{x},\mathbf{x}_{\alpha })=\sum_{n=1}^{N}w_{n}(%
\mathbf{x})\Psi _{n}^{\prime }(\alpha ).  \label{Fourier w alpha}
\end{equation}%
Plugging \eqref{Fourier w} and \eqref{Fourier w alpha} into \eqref{3.4}
gives {\small 
\begin{multline*}
\sum_{n=1}^{N}\partial _{z}w_{n}(\mathbf{x})\Psi _{n}^{\prime }(\alpha )-%
\frac{\partial_{zz}u_{0}(\mathbf{x},\mathbf{x}_{\alpha })}{\partial
_{z}u_{0}(\mathbf{x},\mathbf{x}_{\alpha })}\sum_{n=1}^{N}w_{n}(\mathbf{x}%
)\Psi _{n}^{\prime }(\alpha ) \\
-\frac{\partial }{\partial \alpha }\left( \frac{\partial _{zz}u_{0}(\mathbf{x%
},\mathbf{x}_{\alpha })}{\partial _{z}u_{0}(\mathbf{x},\mathbf{x}_{\alpha })}%
\right) \sum_{n=1}^{N}w_{n}(\mathbf{x})\Psi _{n}(\alpha ) +\sum_{i=1}^{d-1}%
\Big[\frac{\partial _{x_{i}}u_{0}(\mathbf{x},\mathbf{x}_{\alpha })}{\partial
_{z}u_{0}(\mathbf{x},\mathbf{x}_{\alpha })}\sum_{n=1}^{N}\partial
_{x_{i}}w_{n}(\mathbf{x})\Psi _{n}^{\prime }(\alpha ) \\
+ \frac{\partial }{\partial \alpha }\left( \frac{\partial _{x_{i}}u_{0}(%
\mathbf{x},\mathbf{x}_{\alpha })}{\partial _{z}u_{0}(\mathbf{x},\mathbf{x}%
_{\alpha })}\right) \sum_{n=1}^{N}\partial _{x_{i}}w_{n}(\mathbf{x})\Psi
_{n}(\alpha ) -\frac{\partial _{zx_{i}}u_{0}(\mathbf{x},\mathbf{x}_{\alpha
})\partial _{x_{i}}u_{0}(\mathbf{x},\mathbf{x}_{\alpha })}{(\partial
_{z}u_{0}(\mathbf{x},\mathbf{x}_{\alpha }))^{2}}\sum_{n=1}^{N}w_{n}(\mathbf{x%
})\Psi _{n}^{\prime }(\alpha ) \\
-\frac{\partial }{\partial \alpha }\left( \frac{\partial _{zx_{i}}u_{0}(%
\mathbf{x},\mathbf{x}_{\alpha })\partial _{x_{i}}u_{0}(\mathbf{x},\mathbf{x}%
_{\alpha })}{(\partial _{z}u_{0}(\mathbf{x},\mathbf{x}_{\alpha }))^{2}}%
\right) \sum_{n=1}^{N}w_{n}(\mathbf{x})\Psi _{n}(\alpha )\Big]=0.
\end{multline*}%
}

\noindent For each $m\in \{1,\dots ,N\}$, multiply the latter equation by $%
\Psi _{m}(\alpha )$ and then integrate the resulting equation with respect
to $\alpha $. We get 
\begin{equation}
\sum_{n=1}^{N}s_{mn}\partial _{z}w_{n}(\mathbf{x})+\sum_{n=1}^{N}a_{mn}(%
\mathbf{x})w_{n}(\mathbf{x})+\sum_{n=1}^{N}\sum_{i=1}^{d-1}b_{mn,i}(\mathbf{x%
})\partial _{x_{i}}w_{n}(\mathbf{x})=0  \label{3.8}
\end{equation}%
for all $\mathbf{x} \in \Omega$ where $s_{mn}$ is defined as in Proposition %
\ref{prop MK}, 
\begin{multline}
a_{mn}(\mathbf{x})=\int_{-\overline{\alpha }}^{\overline{\alpha }}\left[ -%
\frac{\partial _{zz}u_{0}(\mathbf{x},\mathbf{x}_{\alpha })}{\partial
_{z}u_{0}(\mathbf{x},\mathbf{x}_{\alpha })}\Psi _{n}^{\prime }(\alpha )-%
\frac{\partial }{\partial \alpha }\left( \frac{\partial _{zz}u_{0}(\mathbf{x}%
,\mathbf{x}_{\alpha })}{\partial _{z}u_{0}(\mathbf{x},\mathbf{x}_{\alpha })}%
\right) \Psi _{n}(\alpha )\right. \\
\left. -\sum_{i=1}^{d-1}\frac{\partial }{\partial \alpha }\left( \frac{%
\partial _{zx_{i}}u_{0}(\mathbf{x},\mathbf{x}_{\alpha })\partial
_{x_{i}}u_{0}(\mathbf{x},\mathbf{x}_{\alpha })}{(\partial _{z}u_{0}(\mathbf{x%
},\mathbf{x}_{\alpha }))^{2}}\right) \Psi _{n}^{\prime }(\alpha )\right. \\
\left. -\sum_{i=1}^{d-1}\frac{\partial }{\partial \alpha }\left( \frac{%
\partial _{zx_{i}}u_{0}(\mathbf{x},\mathbf{x}_{\alpha })\partial
_{x_{i}}u_{0}(\mathbf{x},\mathbf{x}_{\alpha })}{(\partial _{z}u_{0}(\mathbf{x%
},\mathbf{x}_{\alpha }))^{2}}\right) \Psi _{n}(\alpha )\right] \Psi
_{m}(\alpha )d\alpha  \label{def matrix A}
\end{multline}%
and for $i=1,\dots, d-1$ 
\begin{equation}
b_{mn,i}(\mathbf{x})=\int_{-\overline \alpha}^{\overline \alpha}\Big[\frac{%
\partial _{x_{i}}u_{0}(\mathbf{x},\mathbf{x}_{\alpha })}{\partial _{z}u_{0}(%
\mathbf{x},\mathbf{x}_{\alpha })}\Psi _{n}^{\prime }(\alpha )+\frac{\partial 
}{\partial \alpha }\left( \frac{\partial _{x_{i}}u_{0}(\mathbf{x},\mathbf{x}%
_{\alpha })}{\partial _{z}u_{0}(\mathbf{x},\mathbf{x}_{\alpha })}\right)
\Psi _{n}(\alpha )\Big]\Psi _{m}(\alpha )d\alpha,  \label{bmn}
\end{equation}%
for all $\mathbf{x}\in \Omega $. For each $\mathbf{x} \in \Omega,$ let $W(%
\mathbf{x})$ $=$$(w_{1}(\mathbf{x}),\dots ,w_{N}(\mathbf{x}))^{T}$, $S$$%
=(s_{mn})_{m,n=1}^{N}$, $A(\mathbf{x})$ $=$$(a_{mn}(\mathbf{x}))_{m,n=1}^{N}$
and $B_{i}(\mathbf{x})=(b_{mn,i}(\mathbf{x}))_{m,n=1}^{N}$ for $i=1,\dots
,d-1$. Since \eqref{3.8} holds true for every $m=1,\dots ,N$, it can be
rewritten as 
\begin{equation}
S_{N}\partial _{z}W(\mathbf{x})+A(\mathbf{x})W(\mathbf{x})+%
\sum_{i=1}^{d-1}B_{i}(\mathbf{x})\partial _{x_{i}}W(\mathbf{x})=0.
\label{No Sinverse}
\end{equation}%
Since $S$ is invertible, see Proposition \ref{prop MK}, then 
\eqref{No
Sinverse} implies the following important system of transport equations 
\begin{equation}
\partial _{z}W(\mathbf{x})+S_{N}^{-1}A(\mathbf{x})W(\mathbf{x}%
)+\sum_{i=1}^{d-1}S_{N}^{-1}B_{i}(\mathbf{x})\partial _{x_{i}}W(\mathbf{x}%
)=0,\quad \mathbf{x}\in \Omega .  \label{3.9}
\end{equation}%
The boundary data for $W$ are: 
\begin{align}
W|_{\partial \Omega} = F(\mathbf{x}) = (f_n)_{n = 1}^N, \quad f_n(\mathbf{x}%
) = \int_{-\overline \alpha}^{\overline \alpha} f(\mathbf{x}, \mathbf{x}%
_\alpha) \partial_z u_0(\mathbf{x}, \mathbf{x}_\alpha) \Psi_n(\alpha)d\alpha
\label{boundary W}
\end{align}
where $f$ is the given data, see \eqref{5}.

\begin{Remark}
From now on, we consider the vector valued function $F(\mathbf{x})$ as the
\textquotedblleft indirect" data, which can be computed directly from %
\eqref{boundary W}. The noiseless data is denoted by $F^*$. The
corresponding noisy data is 
\begin{equation}
F^\delta(\mathbf{x}) = F^*(1 + \delta \mbox{rand}(\mathbf{x})), \quad 
\mathbf{x} \in \partial \Omega  \label{noisy data}
\end{equation}
where $\delta > 0$ is the noise level and rand is a uniformly distributed
function of random numbers taking the range in $[-1, 1].$ \label{rem
indirect data}
\end{Remark}

\begin{Remark}[The approximation context]
Due to the truncation in \eqref{Fourier w}, equation \eqref{3.9} is within
the framework of our approximate mathematical model mentioned in
Introduction. Since this paper is concerned with computational rather than
theoretical results, then this model is acceptable. Our approximation
provides good numerical results in Section \ref{Sec num}. \label{rem 3.2}
\end{Remark}

\begin{Remark}
Problem \ref{problem isp} is reduced to the problem of finding the vector
valued function $W$ satisfying the system \eqref{3.9} and the boundary
condition \eqref{boundary W}. Assume this vector function is computed and
denote it as $W^{\mathrm{comp}}=(w_{1}^{\mathrm{comp}},\dots ,w_{n}^{\mathrm{%
comp}})$. Then, we can compute the function $w^{\mathrm{comp}}(\mathbf{x},%
\mathbf{x}_{\alpha })$ and then the function $u^{\mathrm{comp}}(\mathbf{x},%
\mathbf{x}_{\alpha })$ sequentially via \eqref{Fourier w} and 
\eqref{change
variable}. The computed target function $p^{\mathrm{comp}}(\mathbf{x})$ is
given by \eqref{2.10}.
\end{Remark}

We find an approximate solution of the boundary value problem \eqref{3.9}--%
\eqref{boundary W} by the quasi-reversibility method. This means that we
minimize the functional 
\begin{multline}
J_{\epsilon }(W)=\int_{\Omega }\Big|\partial _{z}W(\mathbf{x}%
)+\sum_{i=1}^{d-1}S_{N}^{-1}B_{i}(\mathbf{x})\partial _{x_{i}}W(\mathbf{x}%
)+S_{N}^{-1}A(\mathbf{x})W(\mathbf{x})\Big|^{2}d\mathbf{x} \\
+\epsilon \Vert W\Vert _{H^{1}(\Omega )^{N}}^{2}  \label{Jepsilon}
\end{multline}%
on the set of vector functions $W\in H^{1}(\Omega )^{N}$ satisfying the
boundary constraint \eqref{boundary W}. Here the space $H^{1}(\Omega )^{N}=%
\underbrace{H^{1}(\Omega )\times \cdot \cdot \cdot \times H^{1}(\Omega )}%
_{N} $ with the commonly defined norm. Similarly to \cite%
{KlibanovAlexeyNguyen:SISC2019}, we analyze the functional $J_{\epsilon }(W)$
for the case when derivatives in \eqref{Jepsilon} are written in finite
differences.

The procedure of computing $p(\mathbf{x})$ is summarized in Algorithm \ref%
{alg}. 
\begin{algorithm}[h!]
\caption{\label{alg}The procedure to solve Problem \ref{problem isp}}
	\begin{algorithmic}[1]
		\State Choose the cut-off number $N = 35$, see Section \ref{Sec num} and Figure \ref{fig choose N}. Find $\{\Psi_n\}_{n = 1}^N.$
		\State Compute the boundary data of the vector valued function $W(\x)$.
		\State \label{Step min}Minimize  the functional $J_{\epsilon}(W)$ subjected to the boundary condition \eqref{boundary W} to obtain $W^{\rm comp}(\x)$, $\x \in \Omega$.
		\State Set $w^{\rm comp}(\x, \x_\alpha) = \sum_{n = 1}^N w_n^{\rm comp} \Psi_n(\alpha)$, $\x \in \Omega$, $\alpha \in [-\overline \alpha, \overline \alpha].$
		\State \label{Step 5}Set $u^{\rm comp} = w^{\rm comp}/\partial_z u_0$. 
		Compute $p^{\rm comp}$ by the average of the left hand side of \eqref{2.10}, namely
		\begin{multline}
		p^{\rm comp} = \frac{1}{2\overline \alpha}\int_{-\overline \alpha}^{\overline \alpha}\Big[ \partial _{z}u_{0}(\mathbf{x},\mathbf{x}_{\alpha })\partial _{z}u^{\rm comp}(\mathbf{x},%
\mathbf{x}_{\alpha })
\\
+\sum_{i=1}^{d-1}\partial _{x_{i}}u_{0}(\mathbf{x},%
\mathbf{x}_{\alpha })\partial _{x_{i}}u^{\rm comp}(\mathbf{x},\mathbf{x}_{\alpha })\Big] d\alpha.
	\label{4.1}
		\end{multline}
	\end{algorithmic}
\end{algorithm}

\section{The quasi-reversibility method in the finite differences}

\label{sec quasi}

For brevity, we describe and analyze here the quasi-reversibility method in
the case when $d=2$. The arguments for higher dimensions can be done in the
same manner. In 2D, $\Omega =(-R,R)\times (a,b).$ We arrange an $N_{x}\times
N_{z}$ grid of points on $\overline{\Omega }$ 
\begin{multline}
\mathcal{G}=\{(x_{i},z_{j}):x_{i}=-R+(i-1)h_{x},z_{j}=a+(j-1)h_{z}, \\
i=1,\dots ,N_{x},j=1,\dots ,N_{z}\},  \label{decompose Omega}
\end{multline}%
where $h_{x}\in \left[ h_{0},\beta _{x}\right) $ and $h_{z}\in \left(
0,\beta _{z}\right) $ are grid step sizes in the $x$ and $z$ directions
respectively and and $h_{0},\beta _{x},\beta _{z}>0$ are certain numbers.
Here, $N_{x}$ and $N_{z}$ are two positive integers. Let $\mathbf{h}=\left(
h_{x},h_{z}\right) .$ We define the discrete set $\Omega ^{\mathbf{h}}$ as
the set of those points of the set \eqref{decompose Omega} which are
interior points of the rectangle $\Omega $ and $\partial \Omega ^{\mathbf{h}%
} $ is the set of those points of the set \eqref{decompose Omega} which are
located on the boundary of $\Omega ,$%
\begin{align*}
\Omega ^{\mathbf{h}}&= \{
(x_{i},z_{j}):x_{i}=-R+(i-1)h_{x},z_{j}=a+(j-1)h_{z}: \\
&\hspace{5cm} i=2,\dots, N_{x}-1;j=2,\dots, N_{z}-1 \} \\
\partial \Omega ^{\mathbf{h}}&=\left\{ (\pm R,z_{j}): j=1,...,N_{z}\right\}
\cup \left\{ (x_{i}, z): i = 1,...,N_{x}, z \in \{a, b\}\right\} , \\
\overline{\Omega }^{\mathbf{h}}&=\Omega ^{\mathbf{h}}\cup \partial \Omega ^{%
\mathbf{h}}.
\end{align*}%
For any continuous function $v$ defined on $\Omega $ its finite difference
version is $v^{\mathbf{h}}=v|_{\mathcal{G}}$. Here, $\mathbf{h}$ denotes the
pair $(h_{x},h_{z}).$ The partial derivatives of the function $v$ are given
via forward finite differences as 
\begin{equation}
\begin{array}{rcl}
\partial _{x}^{h_{x}}v^{\mathbf{h}}(x_{i},z_{j}) & = & \displaystyle\frac{v^{%
\mathbf{h}}(x_{i+1},z_{j})-v^{\mathbf{h}}(x_{i},z_{j})}{h_{x}} \\ 
\partial _{z}^{h_{z}}v(x_{i},z_{j}) & = & \displaystyle \frac{%
v(x_{i},z_{j+1})-v(x_{i},z_{j})}{h_{z}}%
\end{array}
\label{derivative fd}
\end{equation}%
for $i=0,\dots ,N_{x}-1$ and $j=0,\dots ,N_{z}-1.$ We denote the finite
difference analogs of the spaces $L^{2}(\Omega )$ and $H^{1}(\Omega )$ as $%
L^{2,\mathbf{h}}(\Omega )$ and $H^{1,\mathbf{h}}(\Omega )$. Norms in these
spaces are defined as 
\begin{align*}
\Vert v^{\mathbf{h}}\Vert _{L^{2,\mathbf{h}}(\Omega ^{\mathbf{h}})}& =\Big[%
h_{x}h_{z}\sum_{i=0}^{N_{x}}\sum_{j=0}^{N_{z}}\left[ v^{\mathbf{h}%
}(x_{i},z_{j})\right] ^{2}\Big]^{1/2}, \\
\Vert v^{\mathbf{h}}\Vert _{H^{1,\mathbf{h}}(\Omega ^{\mathbf{h}})}& =\Big[%
\Vert v^{\mathbf{h}}\Vert _{L^{2,\mathbf{h}}(\Omega ^{\mathbf{h}%
})}^{2}+h_{x}h_{z}\sum_{i=0}^{N_{x}-1}\sum_{j=0}^{N_{z}-1}\left[ \partial
_{x}^{h_{x}}v^{\mathbf{h}}(x_{i},z_{j})\right] ^{2} \\
& \hspace{6cm}+\left[ \partial _{z}^{h_{z}}v^{\mathbf{h}}(x_{i},z_{j})\right]
^{2}\Big]^{1/2}.
\end{align*}

Let $F^{\mathbf{h}}=F\mid _{\partial \Omega ^{\mathbf{h}}}.$ The problem %
\eqref{3.9}--\eqref{boundary W} becomes 
\begin{multline}
L^{h}\left( W^{\mathbf{h}}\right) = \partial _{z}^{h_{z}}W^{\mathbf{h}%
}(x_{i},z_{j})+S_{N}^{-1}B_{1}(\mathbf{x} _{i},z_{j})\partial _{x}^{h_{x}}W^{%
\mathbf{h}}(x_{i},z_{j}) \\
+S_{N}^{-1}A(x_{i},z_{j})W^{\mathbf{h}}(x_{i},z_{j})=0  \label{51}
\end{multline}
for $i=0,...,N_{x}-1;j=0,...,N_{z} - 1$ and 
\begin{equation}
W^{\mathbf{h}}\mid _{\partial \Omega ^{\mathbf{h}}}=F^{\mathbf{h}}.
\label{53}
\end{equation}%
To solve problem \eqref{51}-\eqref{53} numerically, we introduce the finite
difference version of the functional $J_{\epsilon }$, defined in %
\eqref{Jepsilon}, 
\begin{multline*}
J_{\epsilon }^{\mathbf{h}}(W^{\mathbf{h}})=h_{x}h_{z}\sum_{i=0}^{N_{x}-1}%
\sum_{j=0}^{N_{z}-1}\Big|\partial _{z}^{h_{z}}W^{\mathbf{h}%
}(x_{i},z_{j})+S_{N}^{-1}B_{1}(\mathbf{x}_{i},z_{j})\partial _{x}^{h_{x}}W^{%
\mathbf{h}}(z_{i},z_{j}) \\
+S_{N}^{-1}A(x_{i},z_{j})W^{\mathbf{h}}(x_{i},z_{j})\Big|^{2} +\epsilon
\Vert W^{\mathbf{h}}\Vert _{H_{N}^{1,\mathbf{h}}(\Omega ^{\mathbf{h}})}^{2},
\end{multline*}%
where $H_{N}^{1,\mathbf{h}}(\Omega ^{\mathbf{h}})=\left[ H^{1,\mathbf{h}%
}(\Omega ^{\mathbf{h}})\right] ^{N}$ and similarly for $L_{N}^{2,\mathbf{h}%
}(\Omega ^{\mathbf{h}}).$\ We consider the following problem:

\begin{problem}[Minimization Problem 1]
Minimize the functional $J_{\epsilon }^{\mathbf{h}}(W^{\mathbf{h}})$ on the
set of such vector functions $W^{\mathbf{h}}\in H_{N}^{1,\mathbf{h}}(\Omega
^{\mathbf{h}})$ that satisfy boundary condition \eqref{53}.
\end{problem}

The convergence theory for this problem is formulated in Theorems \ref{thm 1}
and \ref{thm 2}. Proofs of these theorems follow closely the arguments of 
\cite[Section 5]{KlibanovAlexeyNguyen:SISC2019} and are, therefore, not
repeated in this paper. Theorem \ref{thm 1} guarantees the existence and
uniqueness of the minimizer of $J_{\epsilon }^{\mathbf{h}}(W^{\mathbf{h}}),$
and this result can be proven on the basis of Riesz theorem. The next
natural and quite more complicated question is about the convergence of
regularized solutions (i.e. minimizers) to the exact one when the level of
the noise in the data tends to zero, i.e. Theorem \ref{thm 2}. As it is
quite often the case in the quasi-reversibility method (see, e.g. \cite%
{Klibanov:anm2015}), a close analog of Theorem \ref{thm 2} is proven in \cite%
[Section 5]{KlibanovAlexeyNguyen:SISC2019} via applying a new discrete
Carleman estimate: recall that conventional Carleman estimates are in the
continuous form. In other words, these two theorems confirm the
effectiveness of our proposed numerical method for solving Problem \ref%
{problem isp}.

\begin{theorem}[existence and uniqueness of the minimizer]
For any ${\bf h}=(h_{x},h_{z})$
with $h_{x}\in \left[ h_{0},\beta _{x}\right) ,h_{z}\in \left(
0,\beta _{z}\right) ,$any $\epsilon >0$  and for any matrix $%
F^{\mathbf{h}}$ of boundary conditions there exists unique minimizer %
$W_{\min ,\epsilon }^{\mathbf{h}}\in H_{N}^{1,\mathbf{h}}(\Omega ^{\mathbf{h}%
})$ of the functional satisfying boundary condition \eqref{53}.
\label{thm 1}
\end{theorem}

As it is always the case in the regularization theory, assume now that there
exists an \textquotedblleft ideal" solution $W_{\ast }^{\mathbf{h}}\in
H_{N}^{1,\mathbf{h}}(\Omega ^{\mathbf{h}})$ of problem \eqref{51}-\eqref{53}
satisfying the following boundary condition:%
\begin{equation}
W_{\ast }^{\mathbf{h}}\mid _{\partial \Omega ^{\mathbf{h}}}=F_{\ast }^{%
\mathbf{h}},  \label{54}
\end{equation}%
where $F_{\ast }^{\mathbf{h}}$ is the \textquotedblleft ideal" noiseless
boundary data. Since $W_{\ast }^{\mathbf{h}}$ exists, \eqref{54} implies
that there exists an extension $G_{\ast }^{\mathbf{h}}\in H_{N}^{1,\mathbf{h}%
}(\Omega ^{\mathbf{h}})$ with $G_{\ast }^{\mathbf{h}}\mid _{\partial \Omega
^{\mathbf{h}}}=F_{\ast }^{\mathbf{h}}$ of the matrix $F_{\ast }^{\mathbf{h}}$
in $\Omega ^{\mathbf{h}}.$ As to the data $F^{\mathbf{h}}$ in \eqref{53}, we
assume now that there exists an extension $G^{\mathbf{h}}\in H_{N}^{1,%
\mathbf{h}}(\Omega ^{\mathbf{h}})$ with $G^{\mathbf{h}}\mid _{\partial
\Omega ^{\mathbf{h}}}=F^{\mathbf{h}}$ of $F^{\mathbf{h}}$ in $\Omega ^{%
\mathbf{h}}.$ Let $\delta >0$ be the level of the noise in $G^{\mathbf{h}}.$
We assume that 
\begin{equation}
\left\Vert G^{\mathbf{h}}-G_{\ast }^{\mathbf{h}}\right\Vert _{H_{N}^{1,%
\mathbf{h}}(\Omega ^{\mathbf{h}})}<\delta .  \label{55}
\end{equation}%
It is convenient to replace the above notation of the minimizer $W_{\min
,\epsilon }^{\mathbf{h}}$ with $W_{\min ,\epsilon ,\delta }^{\mathbf{h}},$
thus, indicating its dependence on $\delta .$ In \cite[Section 5]%
{KlibanovAlexeyNguyen:SISC2019}, to prove a direct analog of Theorem \ref%
{thm 2} (formulated below), a new Carleman estimate for the finite
difference operator $\partial _{z}^{h_{z}}v$ was proven first. The Carleman
Weight Function of this estimate depends only on the discrete variable $z$.
The value of this function at at the point $z_{j}=a+(j-1)h_{z}$ is $%
e^{2\lambda (j-1)h_{z}},$ where $\lambda >0$ is a parameter. This estimate
is valid only if $\lambda h_{z}<1$ (Lemma 4.7 of \cite[Section 5]%
{KlibanovAlexeyNguyen:SISC2019}). The latter explains the condition of
Theorem \ref{thm 2} imposed on the grid step size $h_{z}$ in the $z-$%
direction.

\begin{theorem}[convergence of regularized solutions]
 Let conditions \eqref{54} and %
\eqref{55} be valid. Let $L^{\mathbf{h}}$ be the operator in 
\eqref{51}. Let $W_{\min ,\epsilon ,\delta }^{\mathbf{h}}$ $\in
H_{N}^{1,\mathbf{h}}(\Omega ^{\mathbf{h}})$ be the minimizer of the
functional $J_{\epsilon }^{\mathbf{h}}(W^{\mathbf{h}})$ with
boundary condition \eqref{53}. Then there exists a sufficiently small number 
$\overline{h}_{z} > 0$ 
depending only on $ h_{0},$ 
$a,$ $b,$ $R,$ $N,$ $L^{\mathbf{h}}$  such that the
following estimate is valid for all $\left( h_{x},h_{z}\right) \in \left[
h_{0},\beta _{x}\right) \times \left( 0,\overline{h}_{z}\right) $ and
all $\epsilon ,\delta >0$ with a constant $C>0$ independent
on $\epsilon ,\delta$%
\begin{equation*}
\left\Vert W_{\min ,\epsilon ,\delta }^{\mathbf{h}}-W_{\ast }^{\mathbf{h}%
}\right\Vert _{L_{N}^{2,\mathbf{h}}\left( \Omega ^{\mathbf{h}}\right) }\leq
C\left( \delta +\sqrt{\epsilon }\left\Vert W_{\ast }^{\mathbf{h}}\right\Vert
_{H_{N}^{1,\mathbf{h}}(\Omega ^{\mathbf{h}})}\right) .
\end{equation*}
\label{thm 2}
\end{theorem}

We also note that Lipschitz stability estimate for problem \eqref{51}-%
\eqref{53} is valid as a direct analog of Theorem 5.5 of \cite[Section 5]%
{KlibanovAlexeyNguyen:SISC2019}. Therefore, uniqueness also takes place for
problem \eqref{51}-\eqref{53}.

\section{Numerical Implementation}

\label{Sec num}

In this section, we solve Problem \ref{problem isp} in the 2D case. The
domain $\Omega $ is 
\begin{equation}
\Omega =(-1,1)\times (1,3).  \label{69}
\end{equation}%
The line of sources $L_{\mathrm{sc}}$ is set to be $(-\overline{\alpha },%
\overline{\alpha })$ with $\overline{\alpha }=3.$

We solve the forward problem to compute the simulated data as follows. Given
the background function $\mathbf{n}_{0}$, instead of solving the nonlinear
Eikonal equation \eqref{eqn Eik0}, we find $u_{0}(\mathbf{x},\mathbf{x}%
_{\alpha })$ using \eqref{def u}. The geodesic line $\Gamma (\mathbf{x},%
\mathbf{x}_{\alpha })$ connecting $\mathbf{x}\in \Omega $ and $\mathbf{x}%
_{\alpha }\in L_{\mathrm{sc}}$ in \eqref{def u} can be found by using the 2D
Fast Marching toolbox which is built in Matlab. The Fast Marching is very
similar to the Dijkstra algorithm to find the shortest paths on graphs. We
refer the reader to \cite{Sethian:cup1999} for more details about Fast
Marching. Next, with this geodesic line in hand, we compute 
\begin{equation*}
u(\mathbf{x},\mathbf{x}_{\alpha })=\int_{\Gamma (\mathbf{x},\mathbf{x}%
_{\alpha })}p(\mathbf{x})d\sigma (\mathbf{x}).
\end{equation*}%
It is clear that the function $u$ solves \eqref{main eqn}. The point $%
\mathbf{x}_{\alpha }$ above is chosen as $(\alpha _{i},0)$ where $\alpha
_{i}=2(i-1)\overline{\alpha }/N_{\alpha }.$ In this paper, we set $N_{\alpha
}=209.$

\begin{figure}[h]
\begin{center}
\subfloat[$N = 10$]{\includegraphics[width=0.3\textwidth]{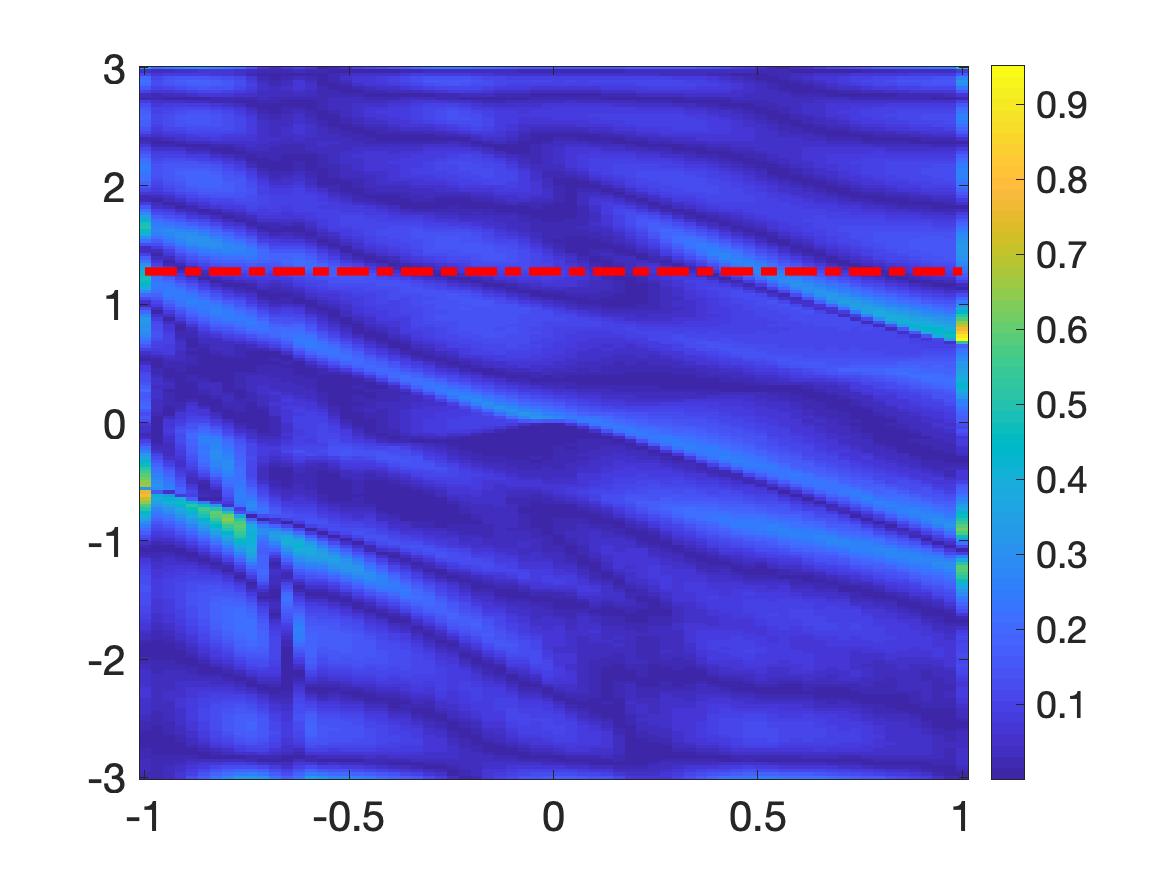}} \hfill %
\subfloat[$N = 20$]{\includegraphics[width=0.3\textwidth]{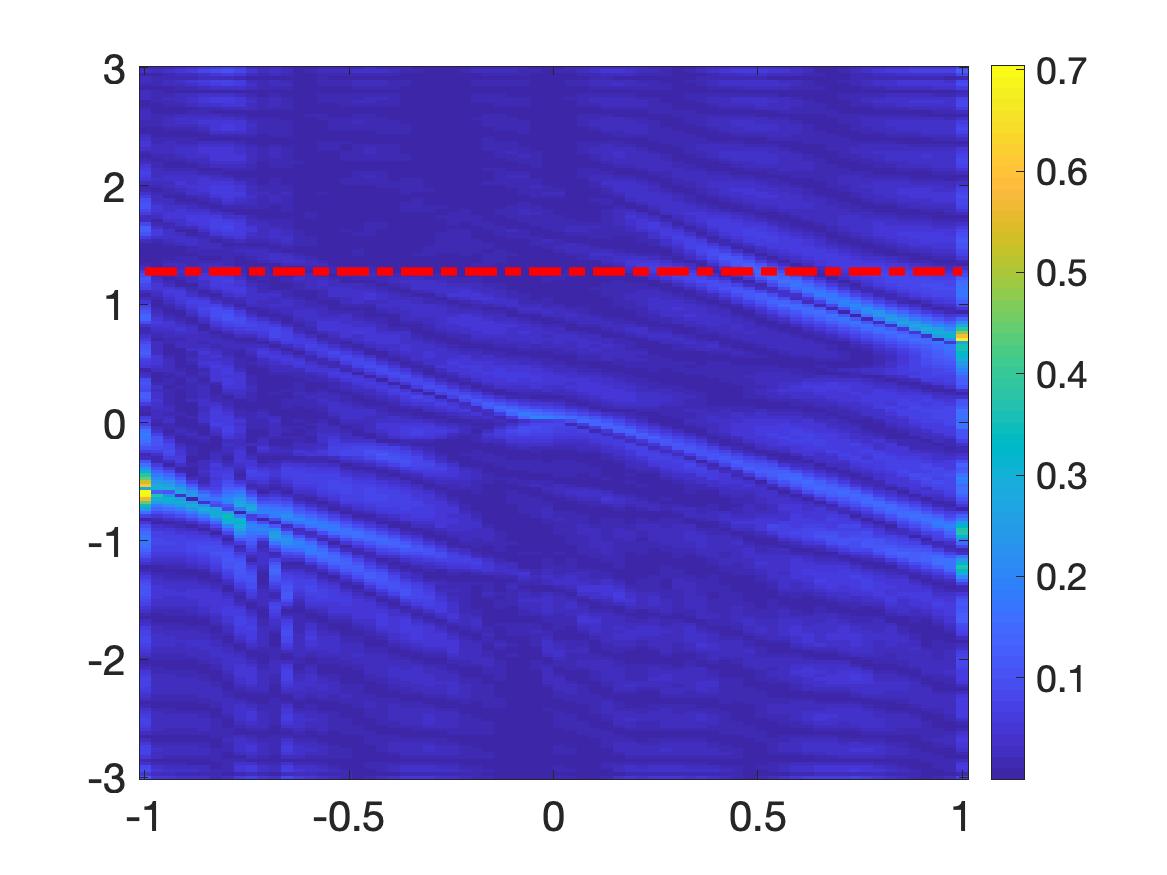}}\hfill %
\subfloat[$N = 35$]{\includegraphics[width=0.3\textwidth]{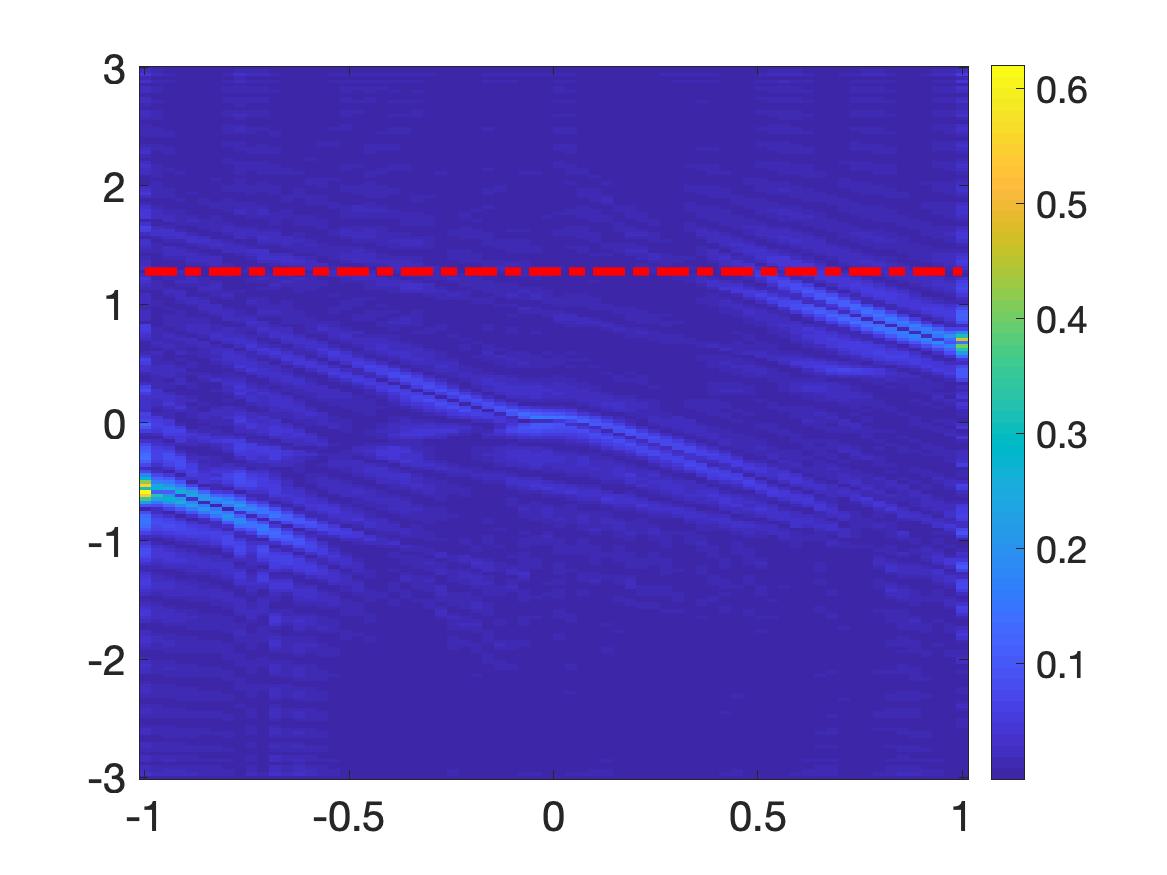}}
\par
\subfloat[$N =
10$]{\includegraphics[width=0.3\textwidth]{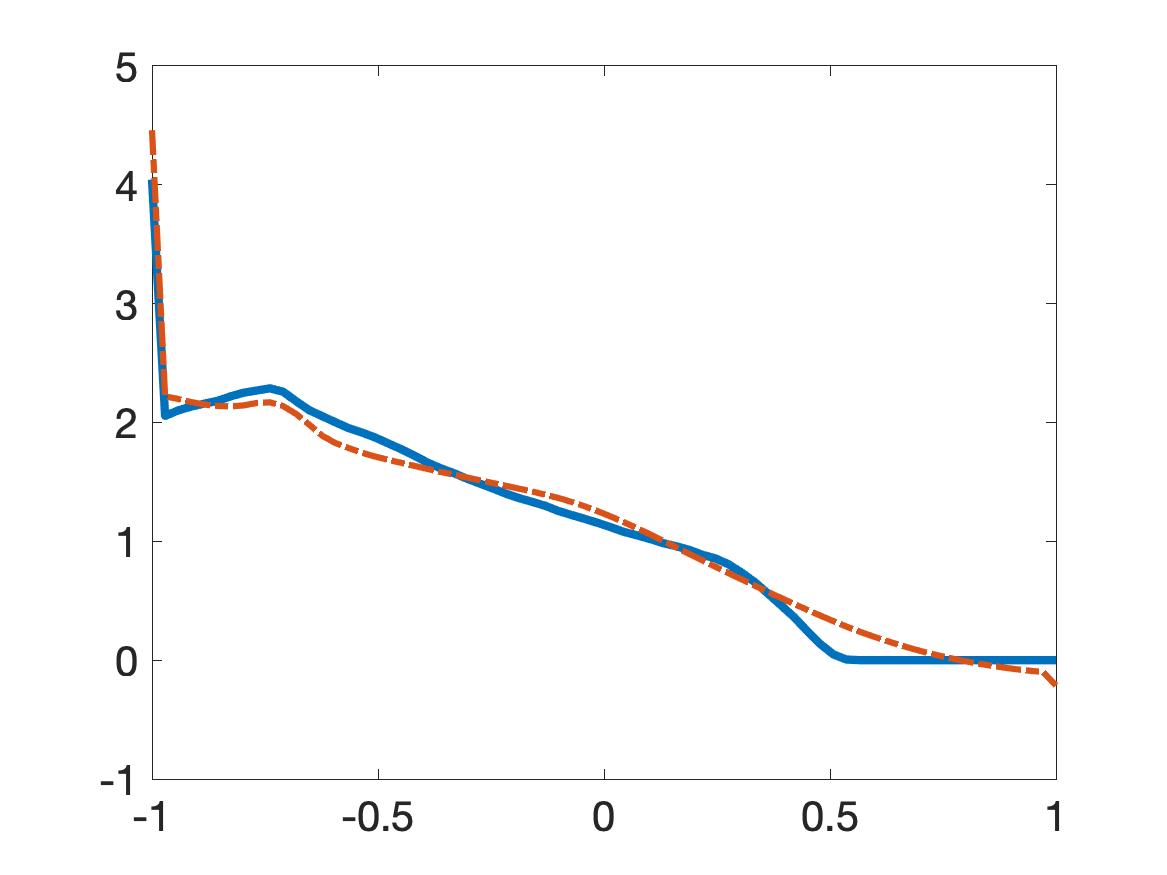}} \hfill 
\subfloat[$N =
20$]{\includegraphics[width=0.3\textwidth]{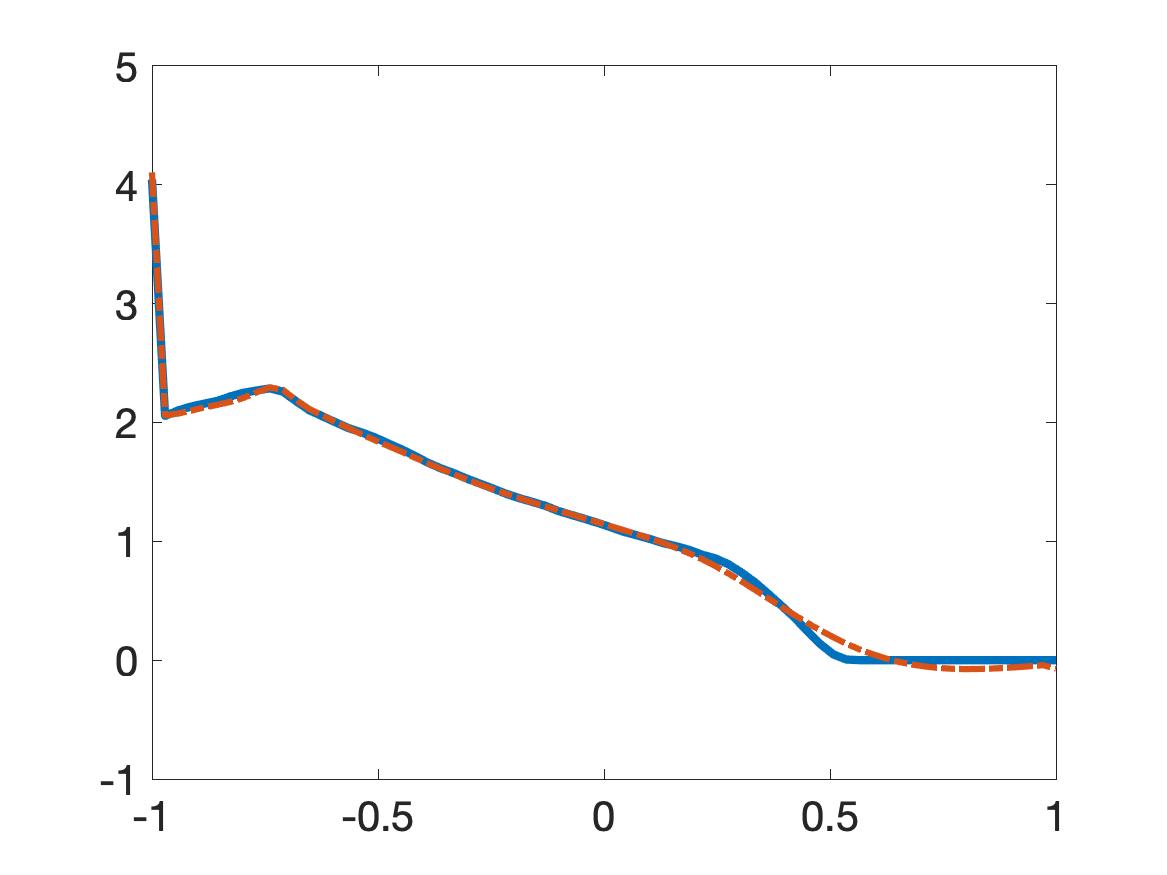}}\hfill 
\subfloat[$N =
35$]{\includegraphics[width=0.3\textwidth]{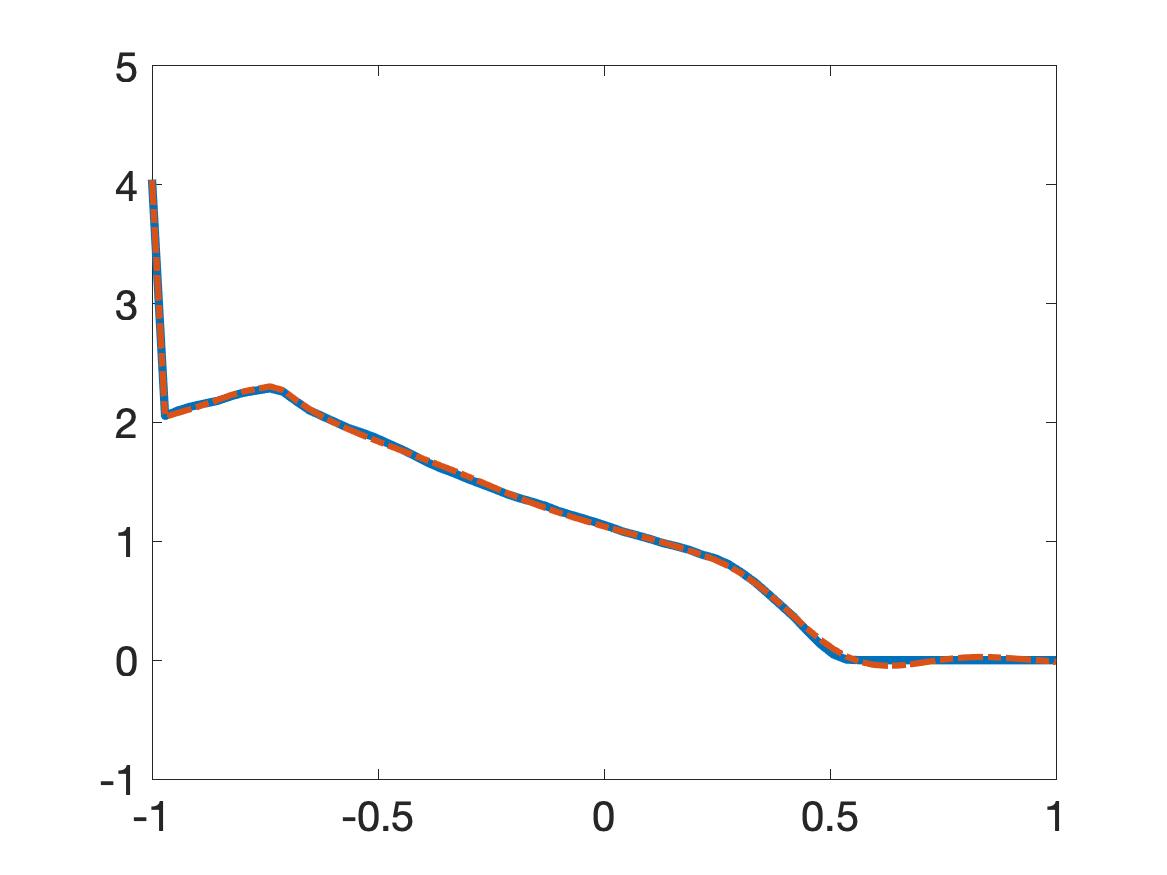}}
\end{center}
\caption{The graphs of the functions $e_{N}(x,\protect\alpha )$ for $N=10,20$
and $35$. The first row shows the 2D graph of $e_{N}$ and the second row
shows the function $w_{\mathrm{true}}(x,z=b,\protect\alpha =1.28)$ and the
function $\sum_{n=1}^{N}w_{n}(x,z=b) \Psi_n(\protect\alpha =1.28)$. We
observe that the larger $N$, the smaller difference of the data and its
approximation is.}
\label{fig choose N}
\end{figure}

We now explain how do we find an appropriate cut-off number $N$ in %
\eqref{Fourier w}. We take the data $f(\mathbf{x},\mathbf{x}_{\alpha })$ in
Test 5 in subsection \ref{sec illu}. %Define the function 
%\[
%	g(x, z = b, \x_\alpha) = w(x, z = b, \x_\alpha) = \partial_z u_0(x, z = b, \x_\alpha) f(x, z = b, \x_\alpha)
%\]
%for $-R \leq x \leq R$, $-\overline \alpha \leq \alpha \leq \overline \alpha$.
%$g = \partial_x u_0(x, z = b, \x_\alpha) f(x, z = b, \x_\alpha)$, which is the function $w$ on $\{(x, z = b, \x_\alpha): -R \leq x \leq R, -\overline \alpha \leq \alpha \leq \overline \alpha\}$.
Then, we compare the function $w(x,z=b,\mathbf{x}_{\alpha })$ and its
approximation $\sum_{n=1}^{N}f_{n}(x,z=b)\Psi _{n}\left( \alpha \right) $
where $f_{n}$ is defined in \eqref{boundary W}. The first row in Figure \ref%
{fig choose N} shows the graphs of 
\begin{equation*}
e_{N}(x,{\alpha })=\Big|w(x,z=b,\mathbf{x}_{{\alpha }})-%
\sum_{n=1}^{N}f_{n}(x,z=b,\mathbf{x}_{{\alpha }})\Big|,\quad x\in
(-R,R),\alpha \in (-\overline{\alpha },\overline{\alpha })
\end{equation*}%
when $N=10,15$ and $35$. The second row in Figure \ref{fig choose N} shows
the true function $w$ and its approximation at $z=b$ and $\alpha =1.28.$ It
is obvious that the sum in the right hand side of the first equation in %
\eqref{Fourier w} when $N=35$ is a good approximation of the data. Thus, we
select $N=35$ in this paper.

\subsection{Computing $W^{\mathrm{comp}}$}

We arrange the grid $\mathcal{G}$ in $\overline{\Omega }$ as in %
\eqref{decompose Omega}. For simplicity, we choose $N_{\mathbf{x}%
}=N_{x}=N_{z}$. The step size $h=h_{x}=h_{z}=2R/(N_{\mathbf{x}}-1)$. We
observe numerically that the matrix $S_{N}^{-1}$, present in the definition
of $J_{\epsilon }$ in \eqref{Jepsilon}, contains some large numbers. This
causes some unwanted errors in computations. Therefore, we slightly modify
the functional $J_{\epsilon }$, see in \eqref{Jepsilon}, suggested by %
\eqref{3.9}, by the following functional (due to \eqref{No Sinverse}) 
\begin{multline}
I_{\epsilon }(W)=\int_{\Omega }|S_{N}\partial _{z}W(\mathbf{x})+A(\mathbf{x}%
)W(\mathbf{x})+\sum_{i=1}^{d-1}B_{i}(\mathbf{x})\partial _{x_{i}}W(\mathbf{x}%
)|^{2}d\mathbf{x} \\
+\epsilon \Vert W\Vert _{H^{1}(\Omega )^{N}}^{2}+\epsilon \Vert \Delta
W\Vert _{L^{2}(\Omega )^{N}}^{2}.  \label{Iepsilon}
\end{multline}%
We have numerically observed that the additional regularization term $%
\epsilon \Vert \Delta W\Vert _{L^{2}(\Omega )^{N}}^{2}$ in \eqref{Iepsilon}
is crucial. Without it, the numerical results do not meet our expectation.
In all tests with all noise level in the data, we choose $\epsilon =10^{-7}$
by a trial and error process. The finite difference version of the
functional $I_{\epsilon }$ for $d=2$ is 
\begin{multline*}
I_{\epsilon }^{h}(W)=h^{2}\sum_{m=1}^{N}\sum_{i,j=1}^{N_{\mathbf{x}}-1}\Big|%
\sum_{n=1}^{N}\Big[\frac{s_{mn}[w_{n}(x_{i},z_{j+1})-w_{n}(x_{i},z_{j})]}{h}
\\
+a_{mn}(x_{i},z_{j})w(x_{i},y_{j})+\frac{%
b_{mn}(x_{i},z_{j})(w(x_{i+1},z_{j})-w(x_{i},z_{j}))}{h}\Big]\Big|^{2} \\
+\epsilon h^{2}\sum_{n=1}^{N}\sum_{i,j=0}^{N_{\mathbf{x}%
}}|w_{n}(x_{i},z_{j})|^{2}+\epsilon h^{2}\sum_{n=1}^{N}\sum_{i,j=0}^{N_{%
\mathbf{x}}-1}\left[ |\partial _{x}^{h}w_{n}(x_{i},z_{j})|^{2}+|\partial
_{z}^{h}w_{n}(x_{i},z_{j})|^{2}\right]
\end{multline*}%
\begin{equation*}
+\epsilon h^{2}\sum_{n=1}^{N}\sum_{i,j=1}^{N_{\mathbf{x}}-1}|\partial
_{xx}w_{n}(x_{i},z_{j})|^{2}+\epsilon h^{2}\sum_{n=1}^{N}\sum_{i,j=1}^{N_{%
\mathbf{x}}-1}|\partial _{zz}^{h}w_{n}(x_{i},z_{j})|^{2}
\end{equation*}%
where $a_{mn}$ and $b_{mn}=b_{mn,1}$ in \eqref{def matrix A} and \eqref{bmn}
respectively. The partial derivatives $\partial _{x}^{h}$ and $\partial
_{z}^{h}$ are as in \eqref{derivative fd}. The second derivatives in finite
difference are understood as usual. %For
%simplicity, we have replaced the quantity $\Vert u_{n}\Vert _{H^{2}(\Omega
%)}^{2}$ in the regularization term by $\Vert u_{n}\Vert _{L^{2}(\Omega
%)}^{2}+\Vert \nabla u_{n}\Vert _{L^{2}(\Omega )}^{2}+\Vert u_{xx}\Vert
%_{L^{2}(\Omega )}+\Vert u_{zz}\Vert _{L^{2}(\Omega )}$ and respectively in
%the finite differences. 
We next line up the discrete vector valued function $w_{n}(x_{i},z_{j})$, $%
1\leq i,j\leq N_{\mathbf{x}}$, $1\leq n\leq N$ as the vector $(\mathfrak{w}_{%
\mathfrak{i}})_{\mathfrak{i}=1}^{N_{\mathbf{x}}^{2}N}$ with 
\begin{equation}
\mathfrak{w}_{\mathfrak{i}}=w_{n}(x_{i},z_{j})  \label{ulineup}
\end{equation}%
where 
\begin{equation}
\mathfrak{i}=(i-1)N_{x}N+(j-1)N+n.  \label{lineup}
\end{equation}%
The functional $I_{\epsilon }^{h}$ in the \textquotedblleft line up" version
is 
\begin{equation}
I_{\epsilon }^{h}(\mathfrak{w})=h^{2}\Big[|\mathcal{L}\mathfrak{w}%
|^{2}+\epsilon |D_{x}\mathfrak{w}|^{2}+\epsilon |D_{y}\mathfrak{u}%
|^{2}+\epsilon |L\mathfrak{w}|^{2}\Big].  \label{6.1}
\end{equation}%
In \eqref{6.1},

\begin{enumerate}
\item $\mathcal{L}$ is the $N_\mathbf{x}^2 N \times N_\mathbf{x}^2 N$ matrix
with entries given by

\begin{enumerate}
\item $(\mathcal{L})_{\mathfrak{i} \mathfrak{j}} = -s_{mn}/h + a_{mn}(x_i,
z_j) - b_{mn}(x_i, y_j)/h$ for $\mathfrak{i} = (i - 1)N_x N + (j - 1)N + m$
and $\mathfrak{j} = (i - 1)N_x N + (j - 1)N + n$,

\item $(\mathcal{L})_{\mathfrak{i} \mathfrak{j}} = b_{mn}(x_i, z_j)/h $ for $%
\mathfrak{i} = (i - 1)N_x N + (j - 1)N + m$ and $\mathfrak{j} = (i + 1 -
1)N_x N + (j - 1)N + n,$

\item $(\mathcal{L})_{\mathfrak{i} \mathfrak{j}} = s_{mn}/h$ for $\mathfrak{i%
} = (i - 1)N_x N + (j - 1)N + m$ and $\mathfrak{j} = (i - 1)N_x N + (j + 1-
1)N + n$,

\item the other entries of $\mathcal{L}$ are $0$
\end{enumerate}

for $1 \leq i, j \leq N_\mathbf{x} - 1$ and $1 \leq m, n \leq N;$

\item $D_x$ is the $N_\mathbf{x}^2 N \times N_\mathbf{x}^2 N$ matrix with
entries given by

\begin{enumerate}
\item $(D_x)_{\mathfrak{i} \mathfrak{i}} = -1/h$ for $\mathfrak{i} = (i -
1)N_x N + (j - 1)N + m$,

\item $(D_x)_{\mathfrak{i} \mathfrak{j}} = 1/h$ for $\mathfrak{i} = (i -
1)N_x N + (j - 1)N + m$ and $\mathfrak{j} = (i + 1 - 1)N_x N + (j - 1)N + m$,

\item the other entries of $\mathcal{L}$ are $0$
\end{enumerate}

for $1 \leq i, j \leq N_\mathbf{x} - 1$ and $1 \leq m \leq N;$

\item $D_y$ is the $N_\mathbf{x}^2 N \times N_\mathbf{x}^2 N$ matrix with
entries given by

\begin{enumerate}
\item $(D_y)_{\mathfrak{i} \mathfrak{i}} = -1/h$ for $\mathfrak{i} = (i -
1)N_x N + (j - 1)N + m$,

\item $(D_y)_{\mathfrak{i} \mathfrak{j}} = 1/h$ for $\mathfrak{i} = (i -
1)N_x N + (j - 1)N + m$ and $\mathfrak{j} = (i - 1)N_x N + (j + 1 - 1)N + m$,

\item the other entries of $\mathcal{L}$ are $0$
\end{enumerate}

for $1 \leq i, j \leq N_\mathbf{x} - 1$ and $1 \leq m \leq N;$

\item $L$ is the $N_\mathbf{x}^2 N \times N_\mathbf{x}^2 N$ matrix with
entries given by

\begin{enumerate}
\item $(L)_{\mathfrak{i} \mathfrak{i}} = -4/h^2$ for $\mathfrak{i} = (i -
1)N_x N + (j - 1)N + m$,

\item $(L)_{\mathfrak{i} \mathfrak{j}} = -1/h^2$ for $\mathfrak{i} = (i -
1)N_x N + (j - 1)N + m$ and $\mathfrak{j} = (i \pm 1 - 1)N_x N + (j - 1)N +
m,$

\item $(L)_{\mathfrak{i} \mathfrak{j}} = -1/h^2$ for $\mathfrak{i} = (i -
1)N_x N + (j - 1)N + m$ and $\mathfrak{j} = (i - 1)N_x N + (j \pm 1 - 1)N +
m,$

\item the other entries of $\mathcal{L}$ are $0$
\end{enumerate}

for $2 \leq i, j \leq N_\mathbf{x} - 1$ and $1 \leq m \leq N.$
\end{enumerate}

The minimizer $\mathfrak{w}$ of $I_{\epsilon }^{h}$ satisfies the equation 
\begin{equation}
\mathcal{L}^{T}\mathcal{L}+\epsilon (\text{Id}%
+D_{x}^{T}D_{x}+D_{y}^{T}D_{y}+L^{T}L)\mathfrak{w}=0.  \label{6.3}
\end{equation}%
On the other hand, due to the constraint \eqref{boundary W} 
\begin{equation}
\mathcal{D}\mathfrak{w}=\mathfrak{f}  \label{6.4}
\end{equation}%
where $\mathcal{D}$ is a $N_{\mathbf{x}}^{2}N\times N_{\mathbf{x}}^{2}N$
matrix and $\mathfrak{f}$ is a $N_{\mathbf{x}}^{2}N$ dimensional vector,
both of which are defined below

\begin{enumerate}
\item $(\mathcal{D})_{\mathfrak{i} \mathfrak{i}} = 1$ for $\mathfrak{i} = (i
- 1)N_x N + (j - 1)N + m$;

\item $(\mathfrak{f})_{\mathfrak{i}} = f_m(x_i, y_j)$ for $\mathfrak{i} = (i
- 1)N_x N + (j - 1)N + m$;

\item the other entries of $\mathcal{L}$ and $\mathfrak{f}$ are $0$
\end{enumerate}

for $i\in \{1,N_{\mathbf{x}}\}$, $1\leq j\leq N_{\mathbf{x}}$ or $2\leq
i\leq N_{\mathbf{x}}-1$, $j\in \{1,N_{\mathbf{x}}\}$ and $1\leq m\leq N.$
Here, $(f_m)_{m = 1}^N$ is in \eqref{boundary W}. Since the data might be
noisy, see \eqref{noisy data}, we slightly modify the system constituted by %
\eqref{6.3} and \eqref{6.4} to a more stable version 
\begin{equation}
\left( \left[ 
\begin{array}{c}
\mathcal{L} \\ 
\mathcal{D}%
\end{array}%
\right] ^{T}\left[ 
\begin{array}{c}
\mathcal{L} \\ 
\mathcal{D}%
\end{array}%
\right] +\epsilon (\text{Id}+D_{x}^{T}D_{x}+D_{y}^{T}D_{y}+L^{T}L)\right) 
\mathfrak{w}=\left[ 
\begin{array}{c}
0 \\ 
\mathfrak{f}%
\end{array}%
\right] .  \label{65}
\end{equation}%
Solving the system \eqref{65}, we obtain $\mathfrak{w}^{\mathrm{comp}}$. The
values of components of vector valued function $W^{\mathrm{comp}}\left( 
\mathbf{x}\right) $ at grid points are computed as $w_{n}(x_{i},z_{j})=%
\mathfrak{w}_{\mathfrak{i}}$ for $\mathfrak{i}=(i-1)N_{x}N+(j-1)N+m$, $1\leq
i,j\leq N_{\mathbf{x}}$, $1\leq m\leq N,$ see \eqref{ulineup}.

We have presented the implementation of Step \ref{Step min} in Algorithm \ref%
{alg}. The other steps are straight forward.

\begin{Remark}[Postprocessing]
In Step \ref{Step 5} of Algorithm \ref{alg} when computing $p^{\mathrm{comp}%
} $ using \eqref{4.1}, which involves $\nabla u^{\mathrm{comp}}$, we smooth
out $u^{\mathrm{comp}}$ by replacing the value of $u^{\mathrm{comp}%
}(x,y,\alpha )$ $\alpha \in \lbrack -\overline{\alpha },\overline{\alpha }]$
by the average of $u^{\mathrm{comp}}$ on the rectangle of $5\times 5$ points
around the point $(x,y)$. We also apply the same smoothing technique for the
function $p^{\mathrm{comp}}.$
\end{Remark}

\subsection{Numerical Tests}

\label{sec illu}

We perform four (4) numerical tests in this paper. When indicating
dependence of any function below on $x,z$, we assume that $\left( x,z\right)
\in \Omega ,$ where the domain $\Omega $ is defined in \eqref{69}\textbf{.}

\begin{Remark}[The function $c_{0}$]
In all our tests below, the function $c_{0}$ is far away from the constant
background function. Therefore, Problem \ref{problem isp} is not considered
as a small perturbation of the problem of inverse Radon transform with
incomplete data, see \cite{KlibanovNguyen:ip2019}. 
All functions $c_{0}$ in our tests
might not smooth in $\R^{2}$ but $c_{0}\in C^{1}(\overline{\Omega })$ in Tests 2,3. Thus, the second derivatives of the corresponding function $u_{0}$ are well-defined in these two tests. 
Even though $c_{0}\notin C^{1}(\overline{\Omega })$ in Test 1, numerically we have not experienced problems with second derivatives of the function $u_{0}$.
\end{Remark}

\begin{figure}[h!]
\begin{center}
\subfloat[The function $p^{\rm
true}$]{\includegraphics[width=0.3\textwidth]{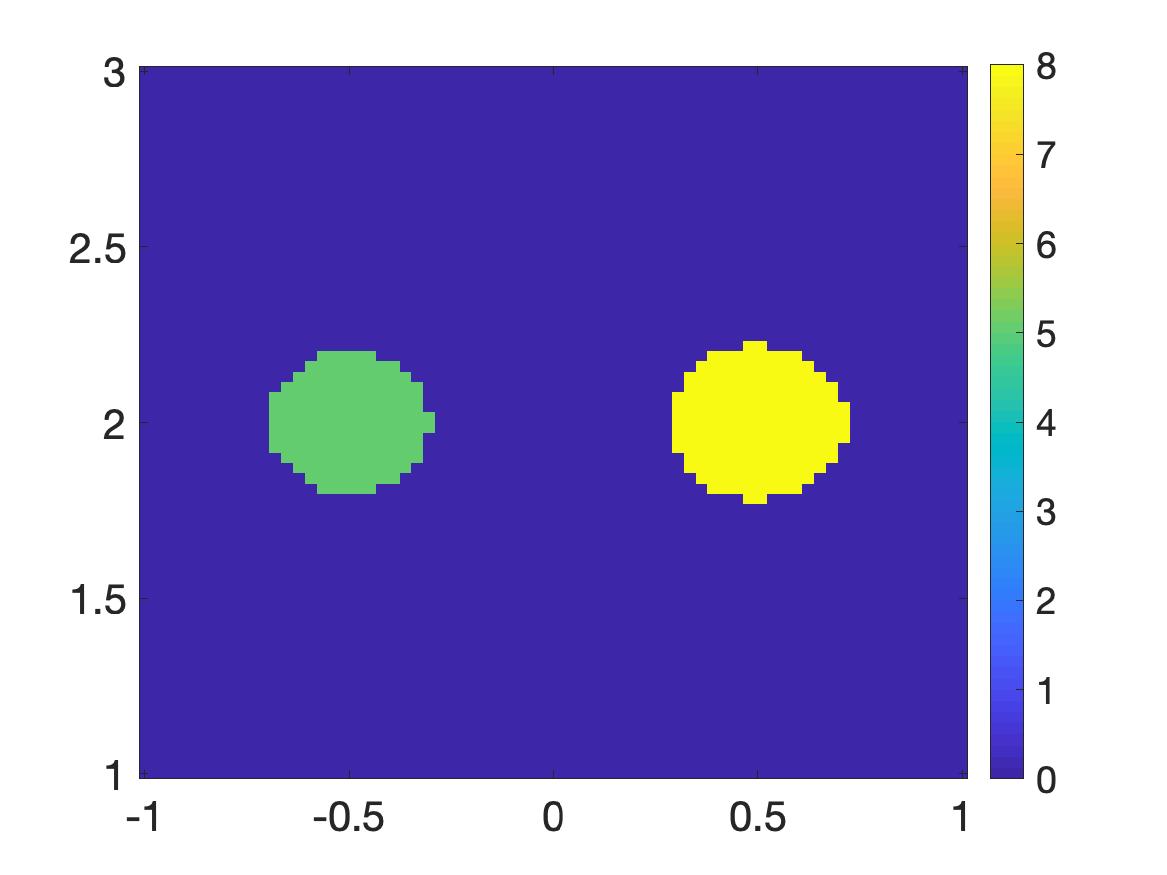}}\hfill 
\subfloat[The function $c_0$ and some geodesic lines, generated by the Fast
Marching package in
Matlab]{\includegraphics[width=0.3\textwidth]{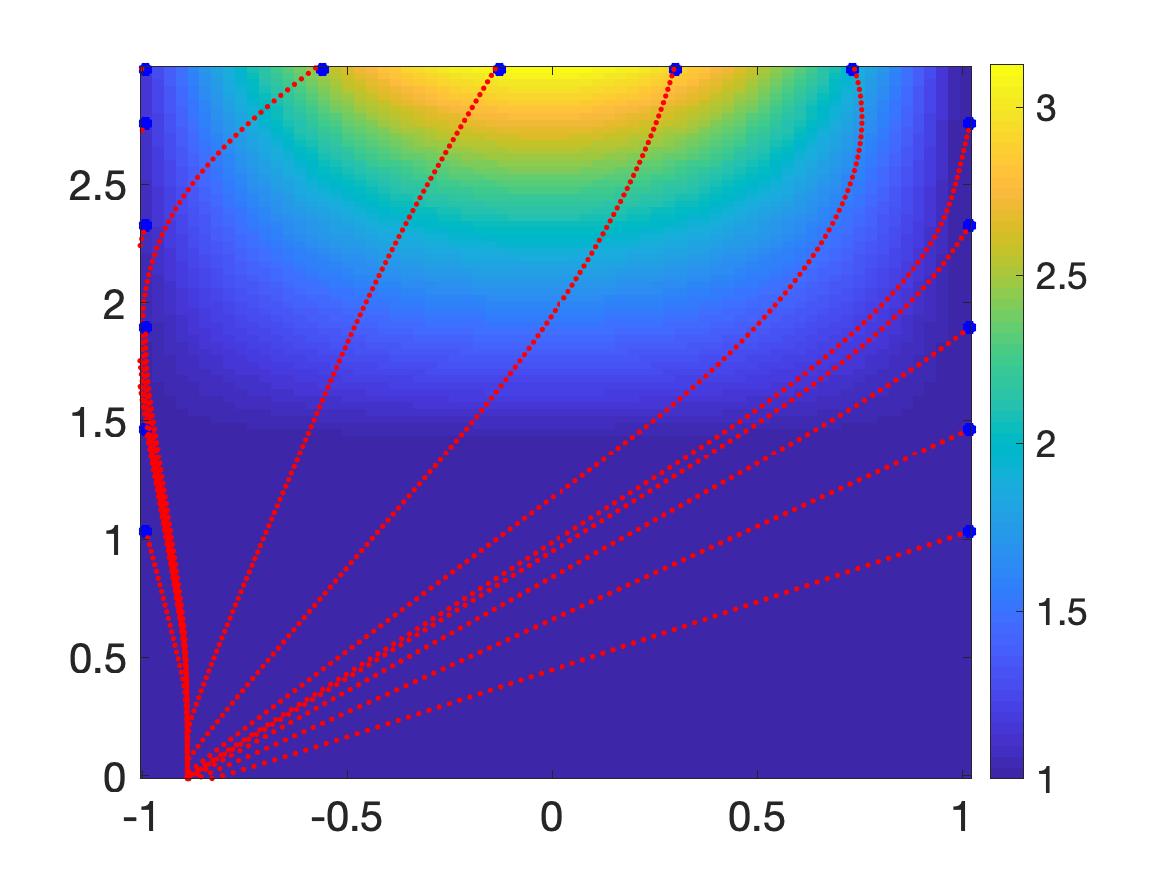}}%
\hfill 
\subfloat[\label{fig test 1 c}The function $p^{\rm comp}$ computed by
Algorithm \ref{alg} with 5\% noise in the data]{\includegraphics[width=0.3\textwidth]{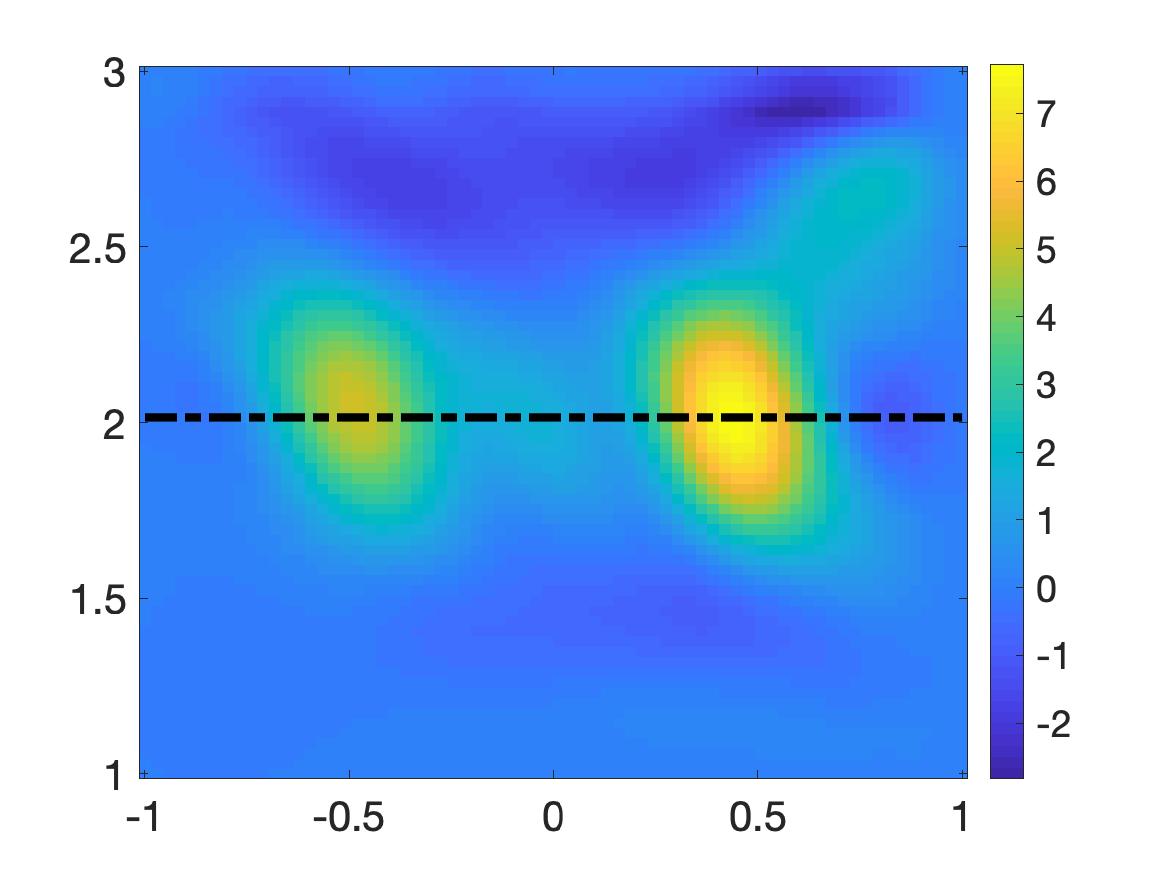}}
\par
\subfloat[The function $p^{\rm true}$ and $p^{\rm comp}$ with 5\% noise in the data
on the set $\{z = 2\},$ indicated by a dash-dot line in
(c)]{\includegraphics[width=0.3\textwidth]{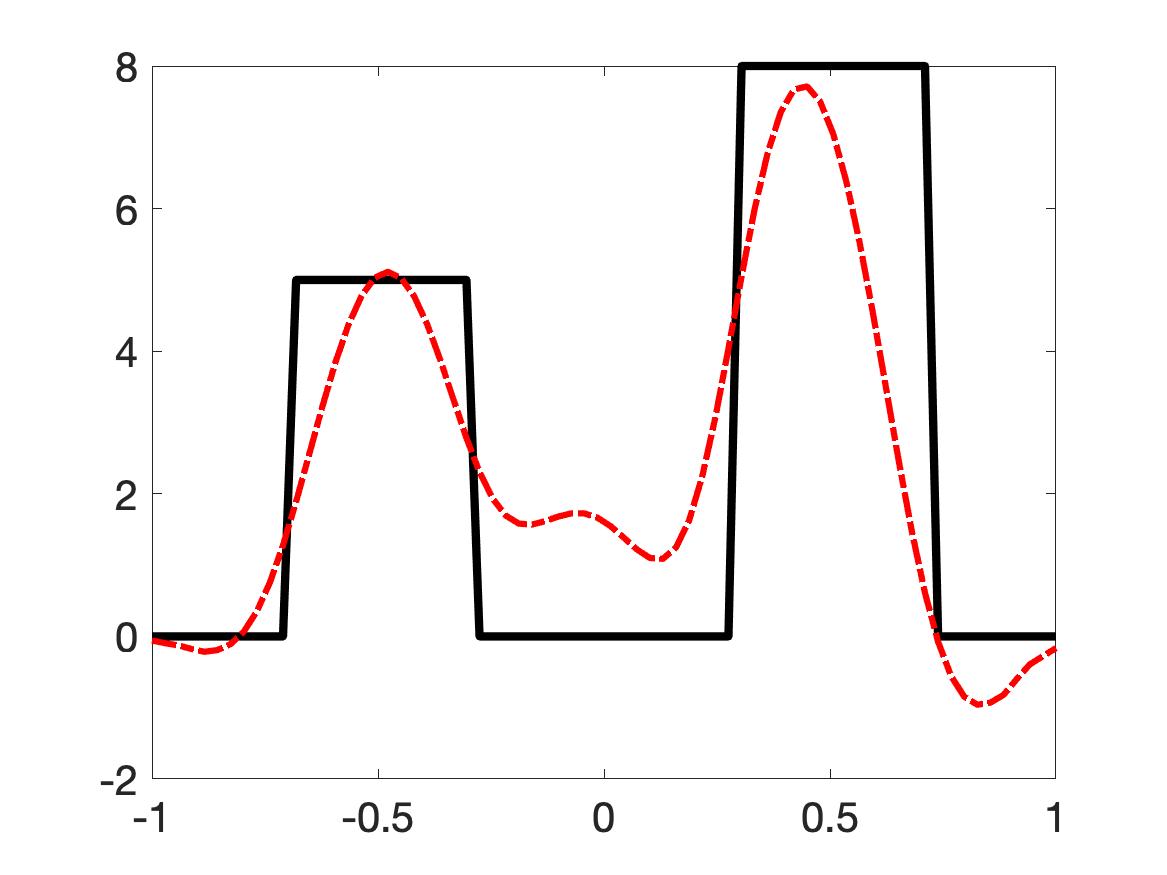}}\hfill 
\subfloat[\label{fig test 1 e}The function $p^{\rm comp}$ computed by Algorithm
\ref{alg} with 120\% noise in the data]{\includegraphics[width=0.3\textwidth]{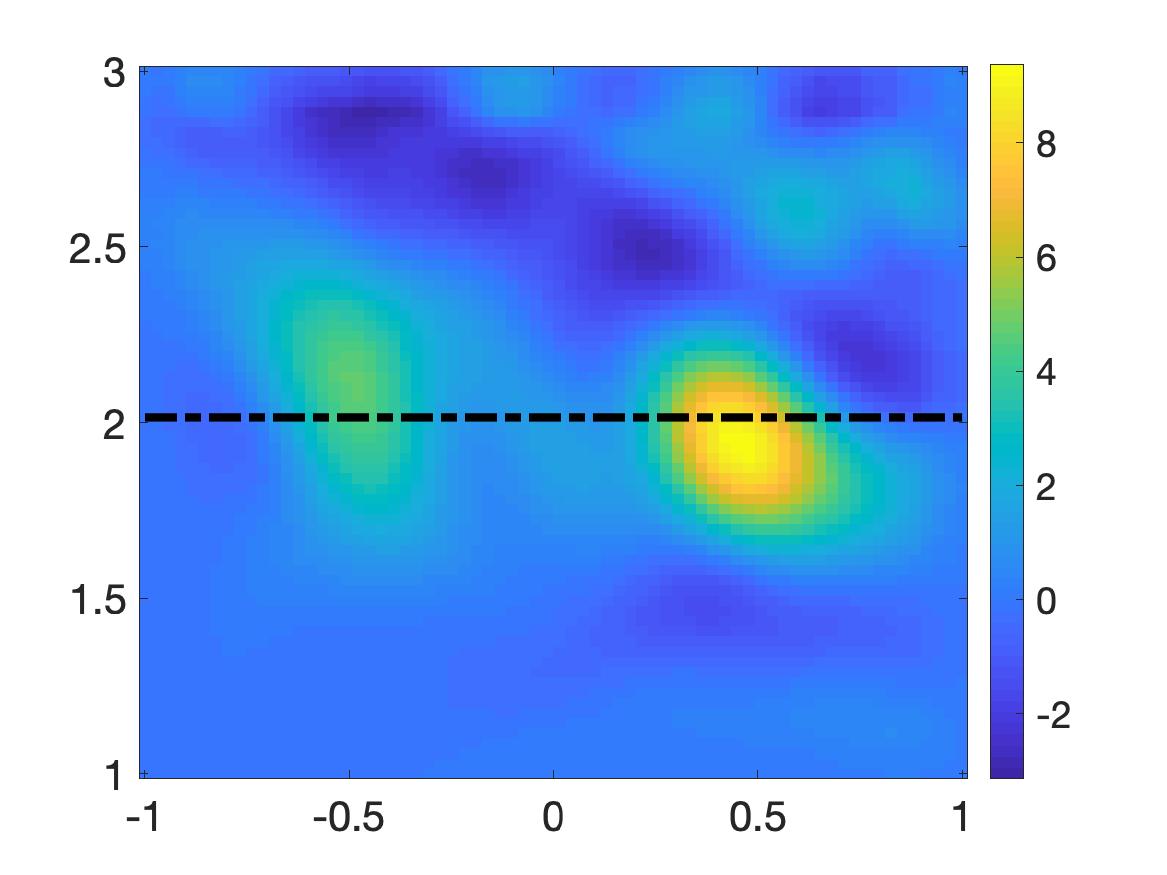}}%
\hfill 
\subfloat[\label{fig test 1 f}The function $p^{\rm true}$ and $p^{\rm comp}$
with 120\% noise in the data on the set $\{z = 2\},$ indicated by a dash-dot line
in (e)]{\includegraphics[width=0.3\textwidth]{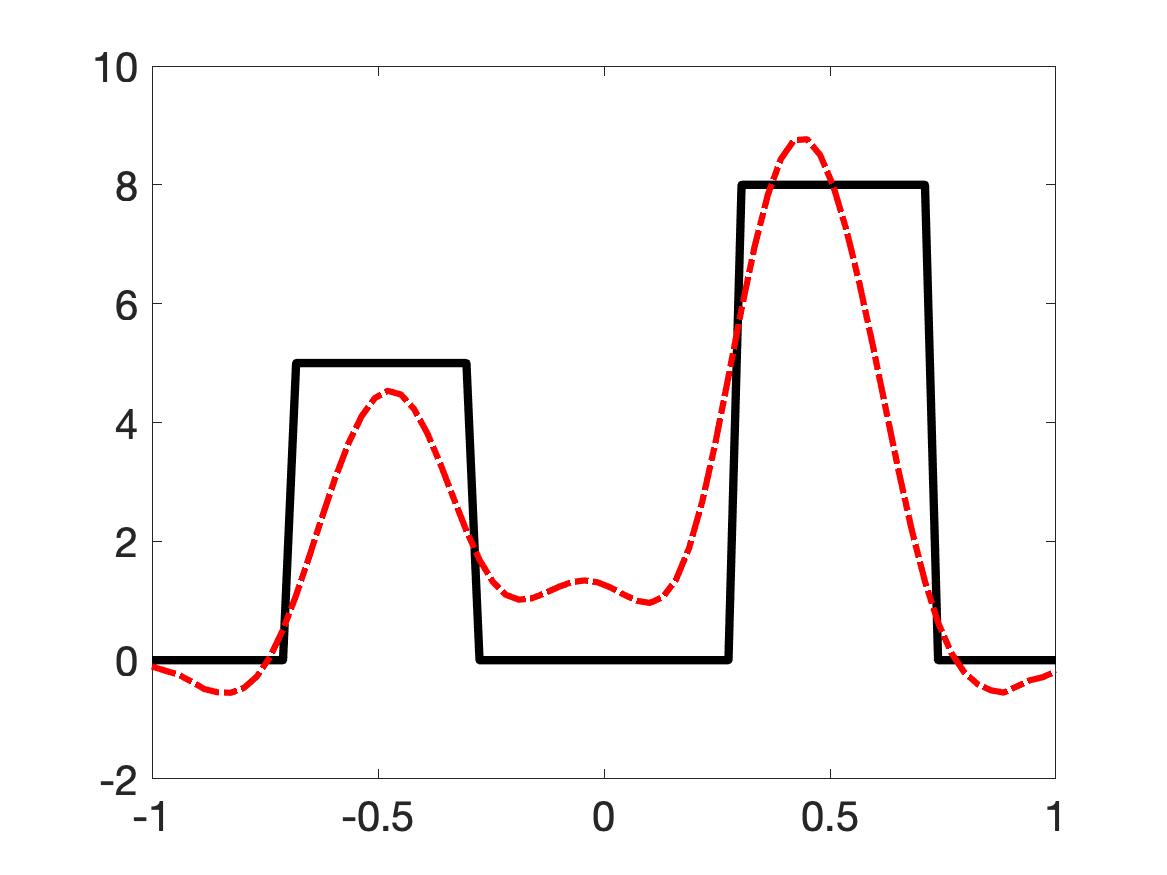}}
\end{center}
\caption{Test 1. The true and reconstructed source functions using Algorithm 
\protect\ref{alg} from noisy data. }
\label{fig test 1}
\end{figure}

\noindent \textbf{Test 1.} The true source function $p$ is given by 
\begin{equation*}
p^{\mathrm{true}}\left( x,z\right) =\left\{ 
\begin{array}{ll}
8 & (x-0.5)^{2}+(z-2)^{2}<0.22^{2}, \\ 
5 & (x+0.5)^{2}+(z-2)^{2}<0.2^{2}, \\ 
0 & \mbox{otherwise.}%
\end{array}%
\right.
\end{equation*}%
The background function $c_{0}$ is 
\begin{equation*}
c_{0}\left( x,z\right) =\left\{ 
\begin{array}{ll}
1+0.3(1-x^{2})(z^{2}-2) & \text{if }z^{2}-2>0, \\ 
1 & \mbox{otherwise}.%
\end{array}%
\right.
\end{equation*}%
The numerical results of this test are displayed in Figure \ref{fig test 1}.

The support of $p^{\mathrm{true}}$ in Test 1 consists of two discs. The
value of the function $p$ in the right disc is higher than the value in the
left disc. Our method detects both these inclusions very well, see Figures %
\ref{fig test 1 c}--\ref{fig test 1 f}. There are some unwanted artifacts
near $\partial \Omega$ where we measure the noisy data. The higher level of
noisy data, the more artifacts present. When the noise level $\delta =5\%$,
the computed maximal value of $p^{\text{comp}}$ in the left inclusion is
5.16 (relative error 3.2\%) and the computed maximal value of $p^{\text{comp}%
}$ in the right inclusion is 7.72 (relative error 3.5\%). When the noise
level $\delta =120\%$, the computed maximal value of $p^{\text{comp}}$ in
the left inclusion is 4.71 (relative error 5.8\%) and the computed maximal
value of $p^{\text{comp}}$ in the right inclusion is 9.37 (relative error
17.1\%).

\begin{figure}[h!]
\begin{center}
\subfloat[The function $p_{\rm
true}$]{\includegraphics[width=0.3\textwidth]{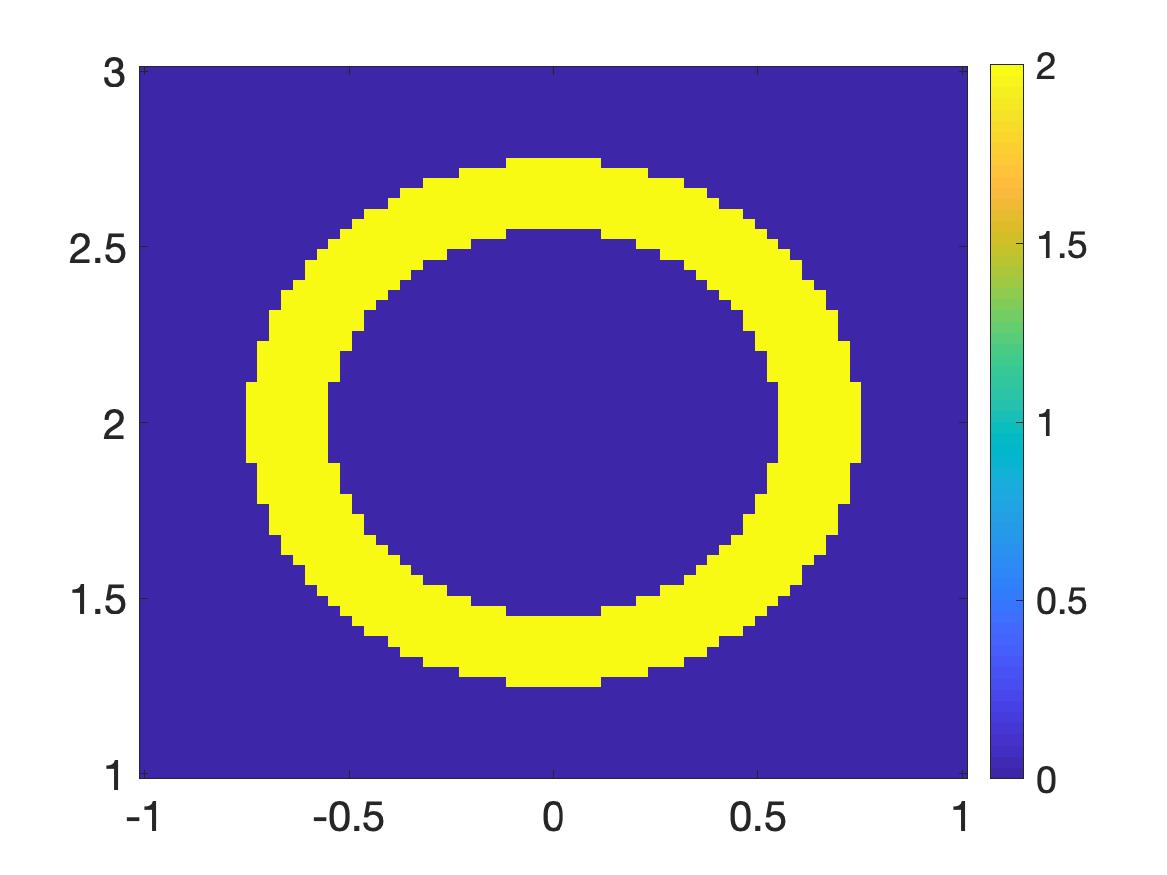}}\hfill 
\subfloat[The function $c_0$ and some geodesic line, generated by the Fast
Marching package in
Matlab]{\includegraphics[width=0.3\textwidth]{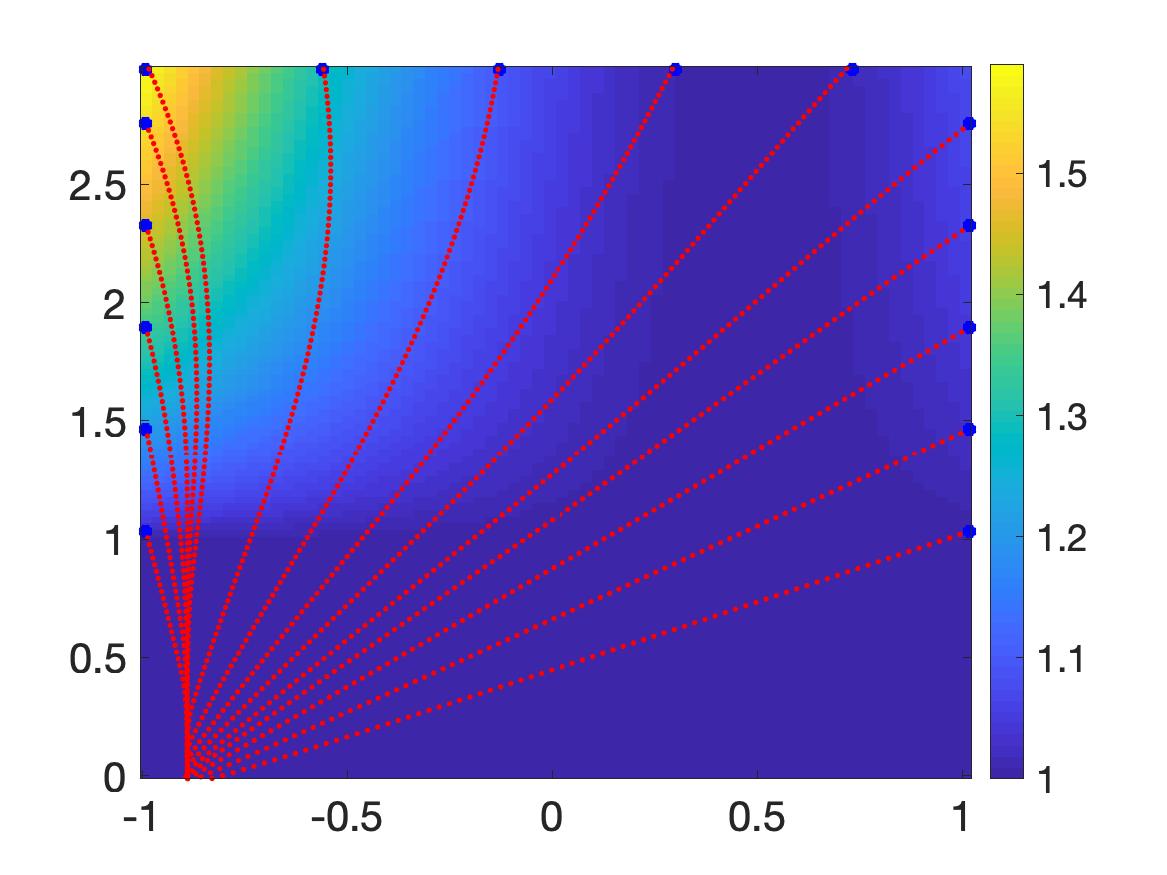}}%
\hfill 
\subfloat[\label{fig test 3 c}The function $p^{\rm comp}$ computed by
Algorithm \ref{alg} with 5\% noise in the
data]{\includegraphics[width=0.3\textwidth]{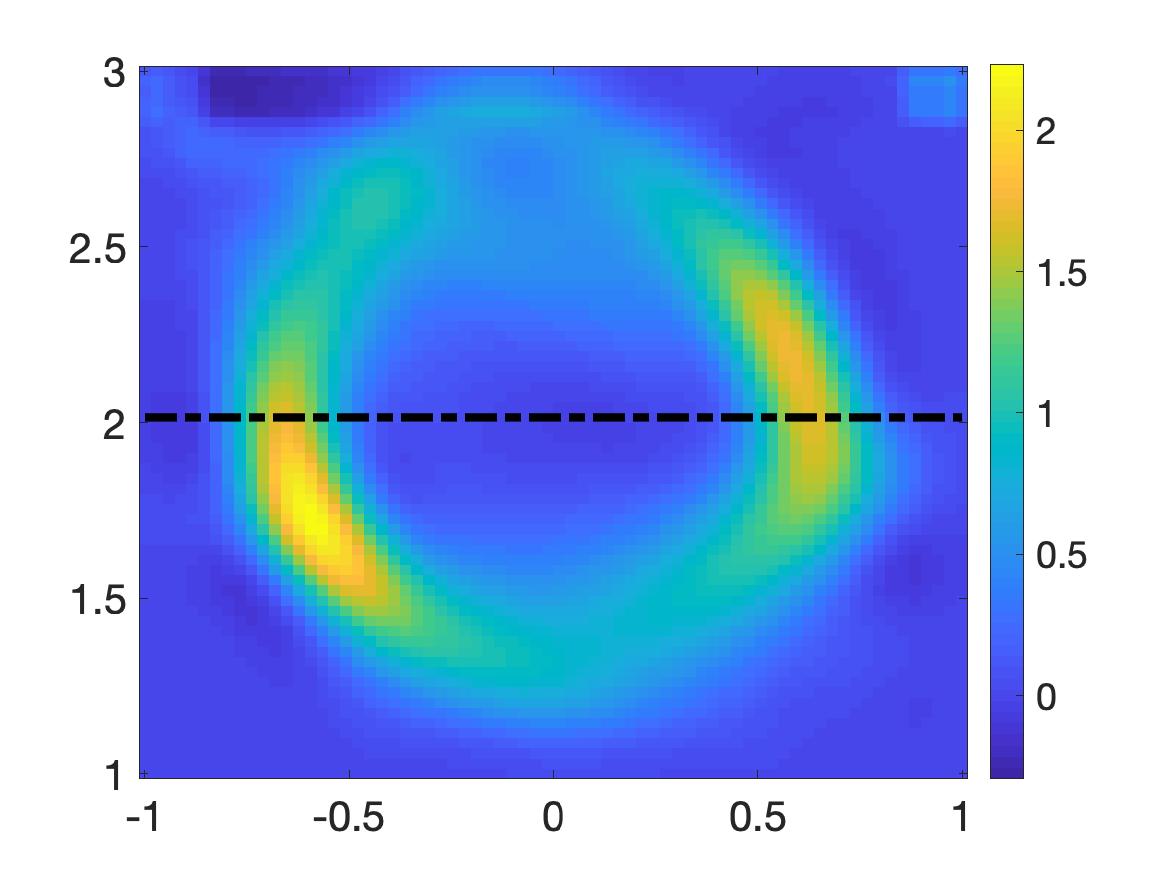}}
\par
\subfloat[\label{fig test 3 d}The function $p^{\rm true}$ and $p^{\rm comp}$
with 5\% noise in the data on the set $\{z = 2\},$ indicated by a dash-dot line in
(c)]{\includegraphics[width=0.3\textwidth]{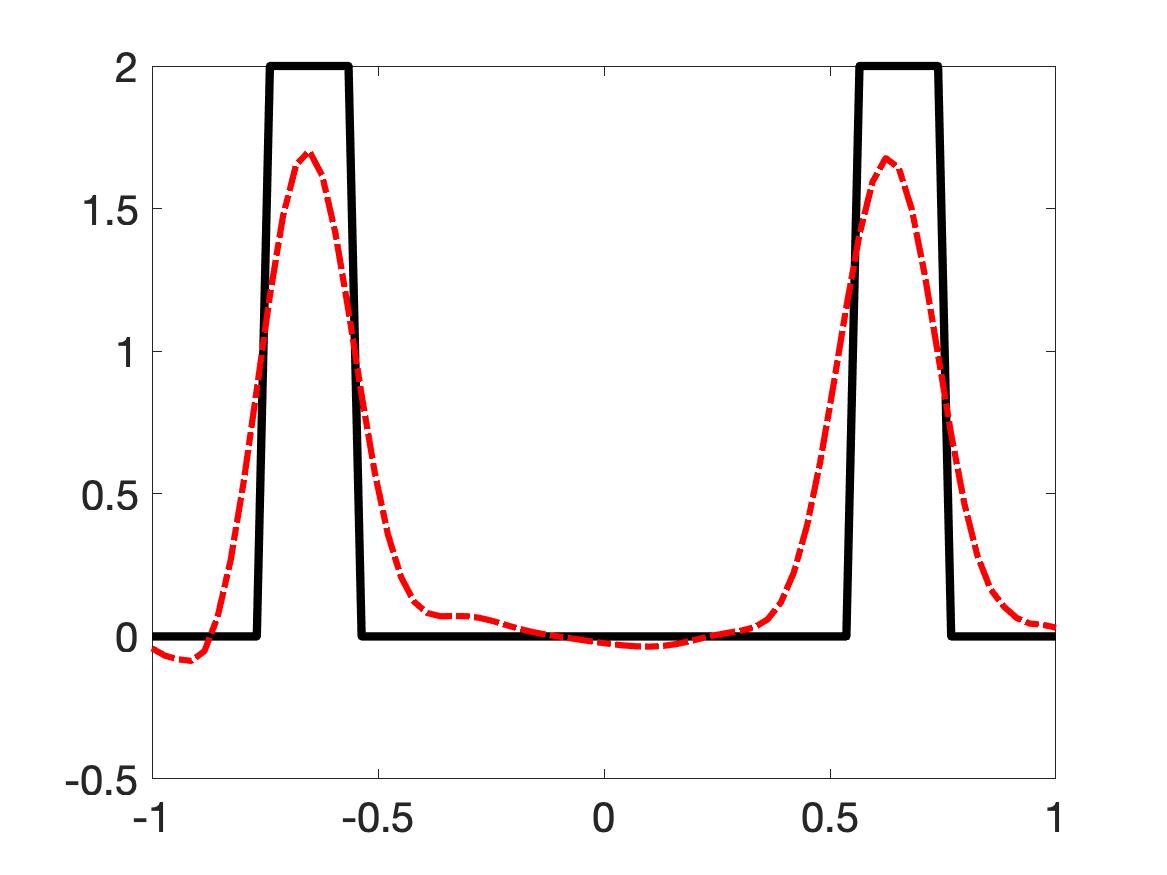}}\hfill 
\subfloat[\label{fig test 3 e}The function $p^{\rm comp}$ computed by Algorithm
\ref{alg} with 30\% noise in the
data]{\includegraphics[width=0.3\textwidth]{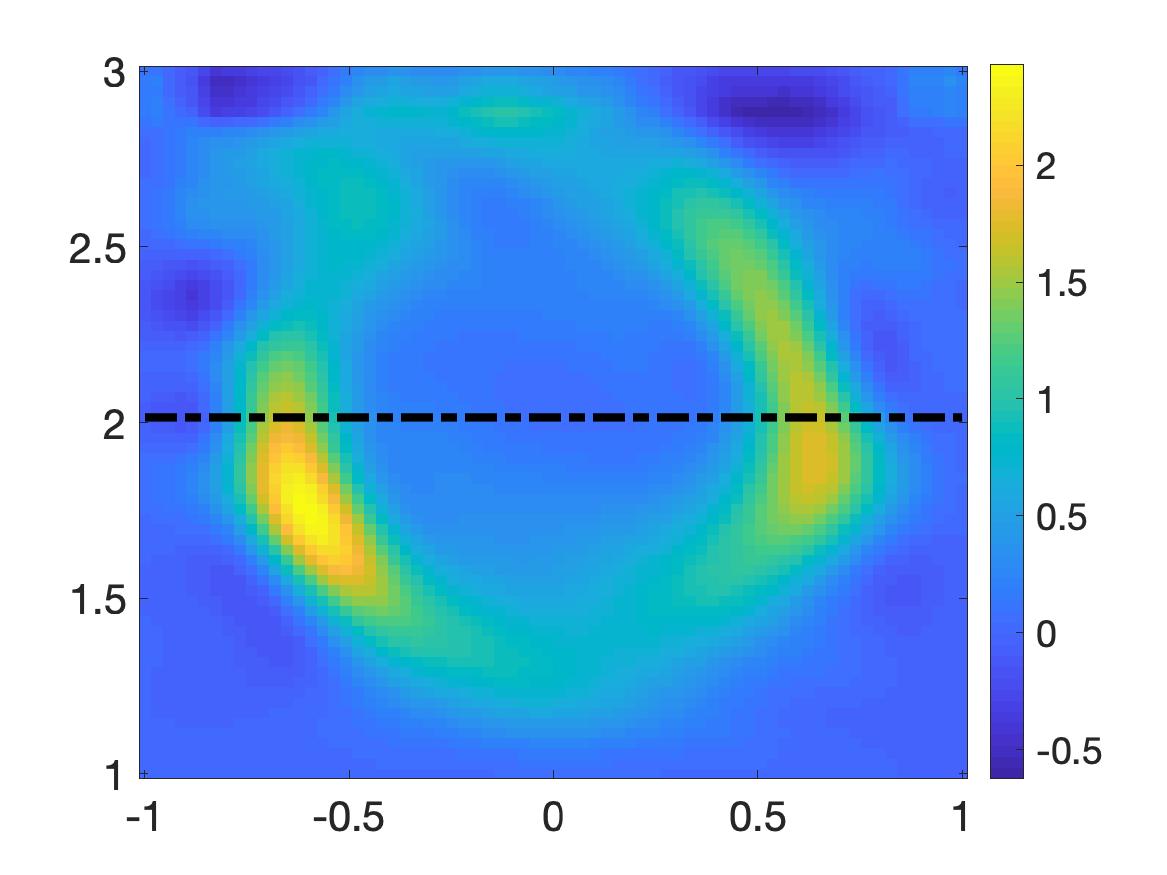}}\hfill 
\subfloat[\label{fig test 3 f}The function $p^{\rm true}$ and $p^{\rm comp}$
with 30\% noise in the data on the set $\{z = 2\},$ indicated by a dash-dot line in
(e)]{\includegraphics[width=0.3\textwidth]{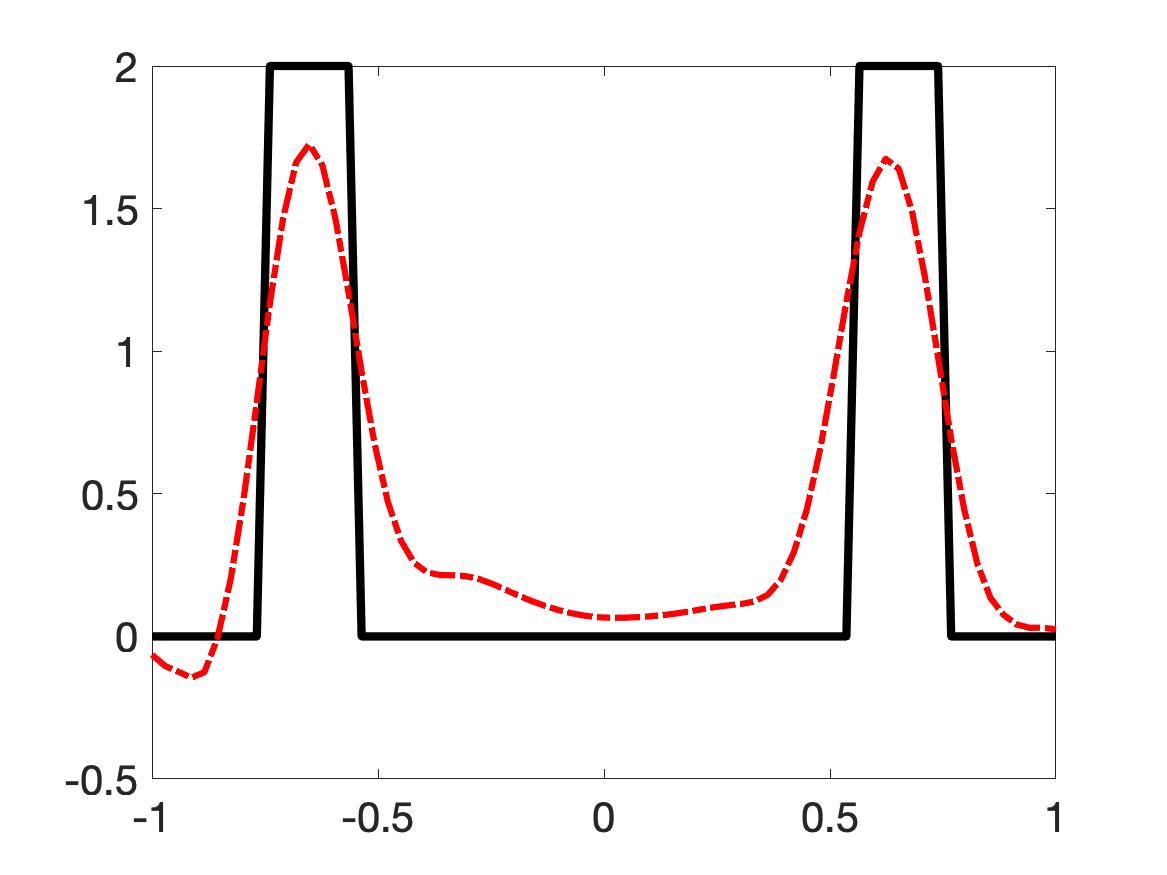}}
\end{center}
\caption{Test 2. The true and reconstructed source functions using Algorithm 
\protect\ref{alg} from noisy data. }
\label{fig test 3}
\end{figure}

\noindent \textbf{Test 2.} We test a complicated case when the support of $%
p_{\mathrm{true}}$ looks like a ring. In this test, 
\begin{equation*}
p^{\mathrm{true}}\left( x,z\right) =\left\{ 
\begin{array}{ll}
2 & 0.55^{2}<r^{2}=x^{2}+(z-2)^{2}<0.75^{2}, \\ 
0 & \mbox{otherwise.}%
\end{array}%
\right.
\end{equation*}%
The background function $c_{0}$ is given by 
\begin{equation*}
c_{0}\left( x,z\right) =\left\{ 
\begin{array}{ll}
1+0.25(x-0.5)^{2}\ln (z) & z>1, \\ 
1 & \mbox{otherwise.}%
\end{array}%
\right.
\end{equation*}%
The numerical results of this test are displayed in Figure \ref{fig test 3}.

In this test, it is evident that the reconstructed \textquotedblleft ring"
is acceptable, see Figures \ref{fig test 3 c} and \ref{fig test 3 e}. The
position of the ring is detected quite well, see Figures \ref{fig test 3 d}
and \ref{fig test 3 f}. When the noise level is $5\%$, the reconstructed
maximal value of $p^{\text{comp}}$ in the ring is 2.23 (relative error
11.5\%). When the noise level is $30\%$, the reconstructed maximal value of $%
p^{\text{comp}}$ in the ring is 2.42 (relative error 21.0\%).

\begin{figure}[h!]
\begin{center}
\subfloat[The function $p^{\rm
true}$]{\includegraphics[width=0.3\textwidth]{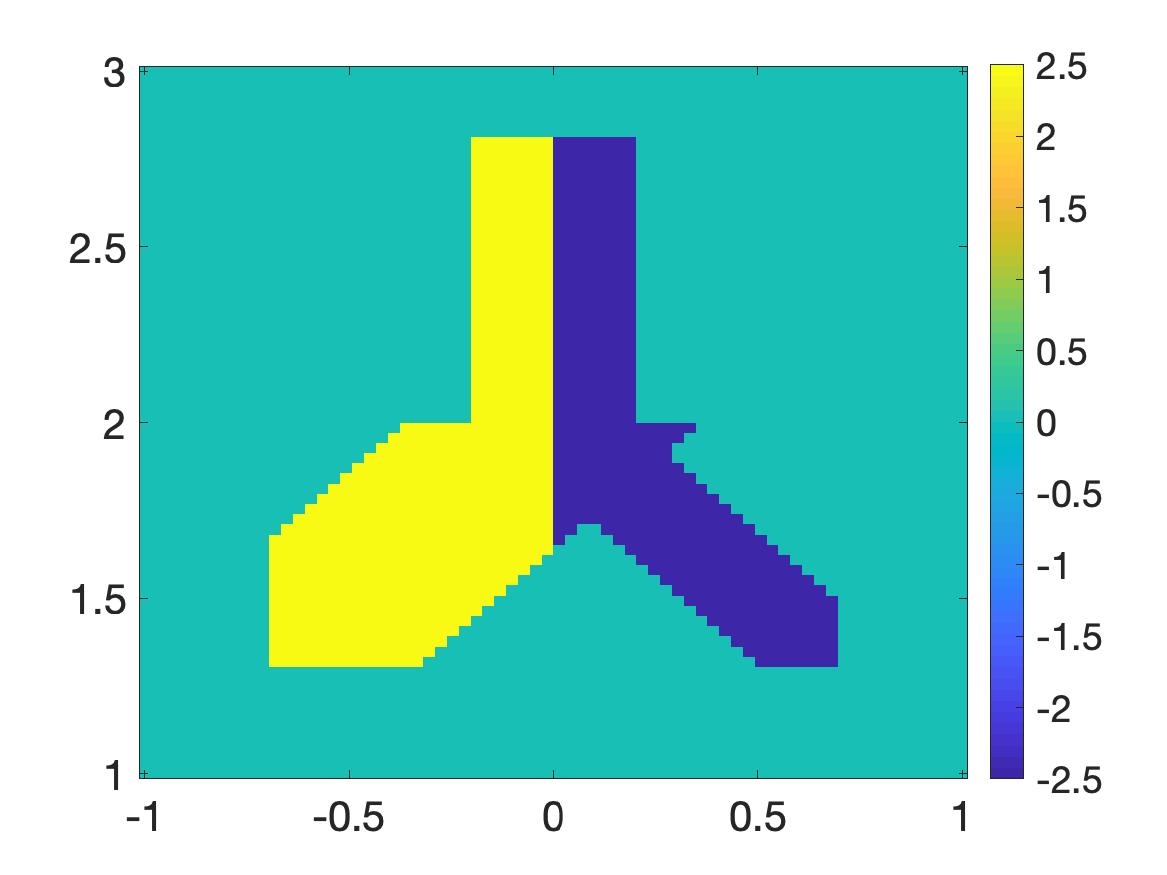}}\hfill 
\subfloat[The function $c_0$ and some geodesic line, generated by the Fast
Marching package in
Matlab]{\includegraphics[width=0.3\textwidth]{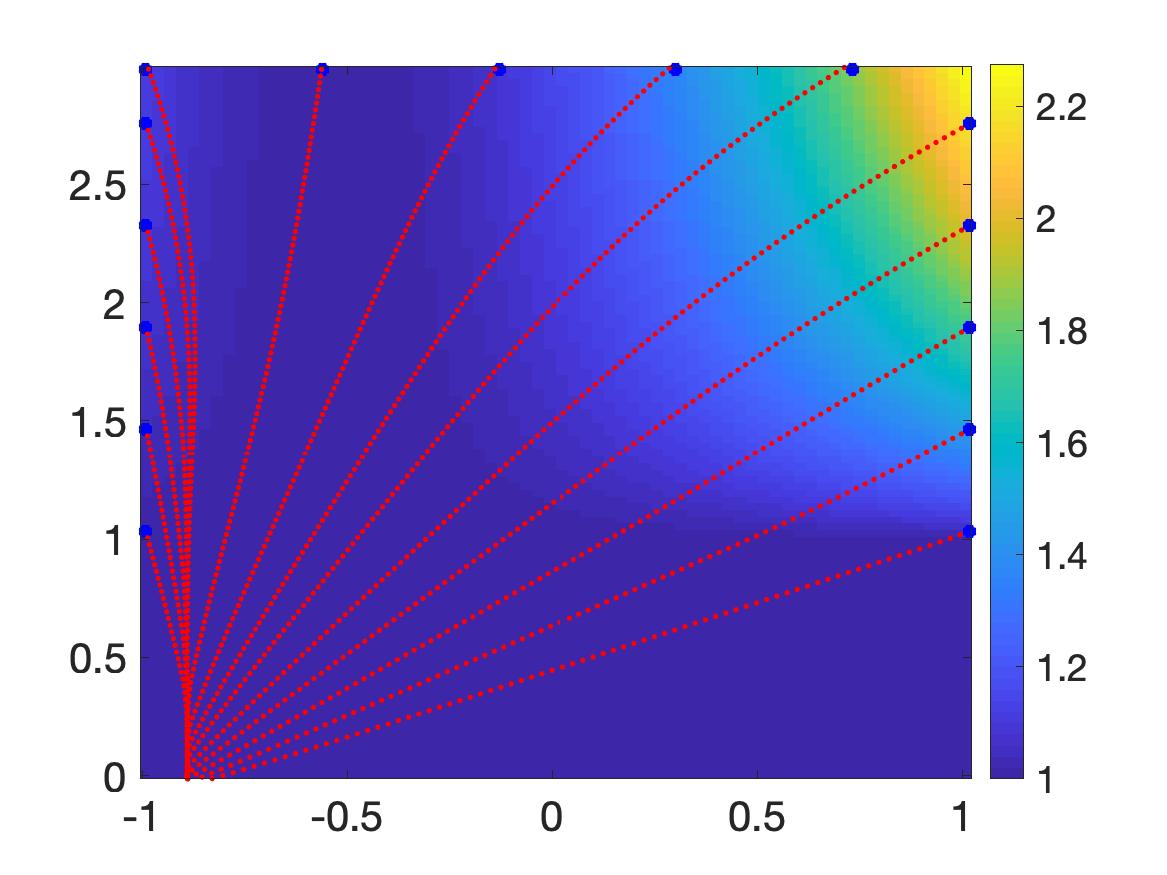}}%
\hfill 
\subfloat[\label{fig test 4 c}The function $p^{\rm comp}$ computed by
Algorithm \ref{alg} with 5\% noise in the
data]{\includegraphics[width=0.3\textwidth]{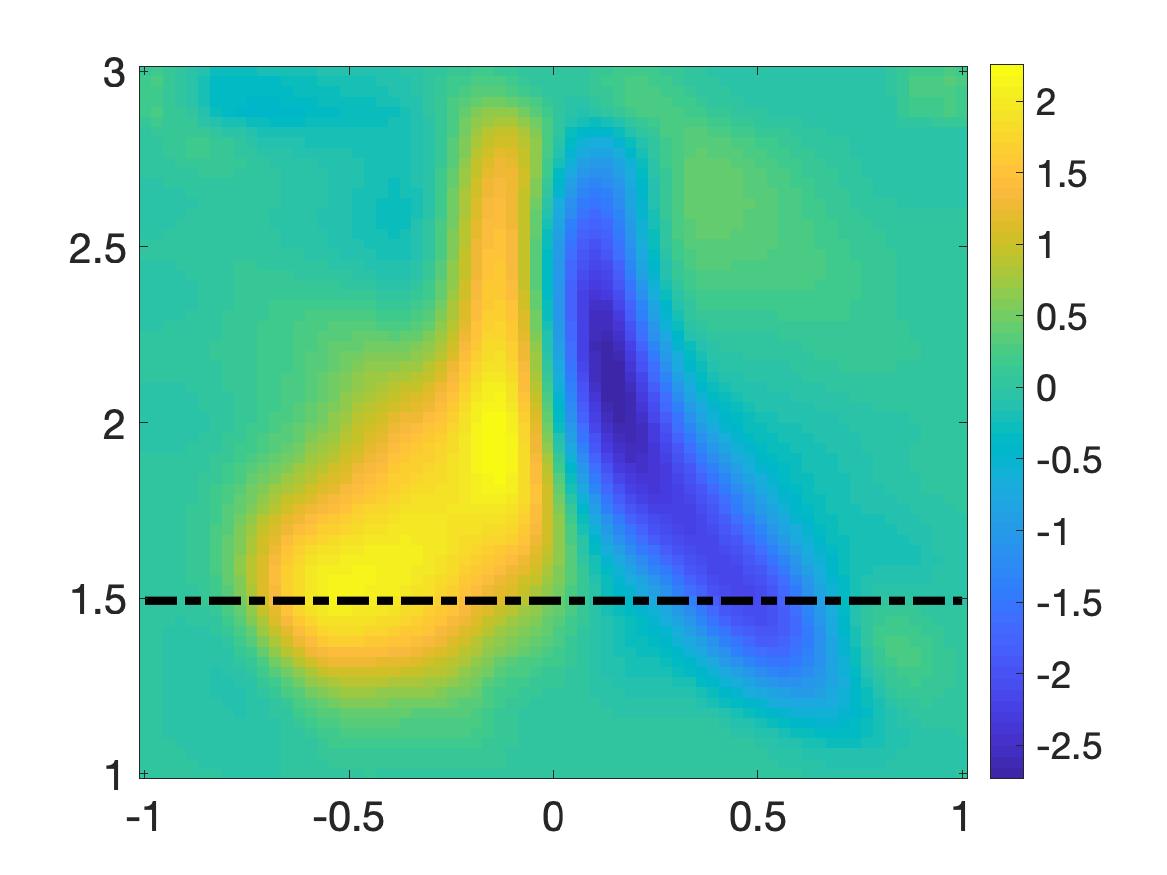}}
\par
\subfloat[\label{fig test 4 d}The function $p^{\rm true}$ and $p^{\rm comp}$
with 5\% noise in the data on the set $\{z = 1.5\},$ indicated by a dash-dot line
in (c)]{\includegraphics[width=0.3\textwidth]{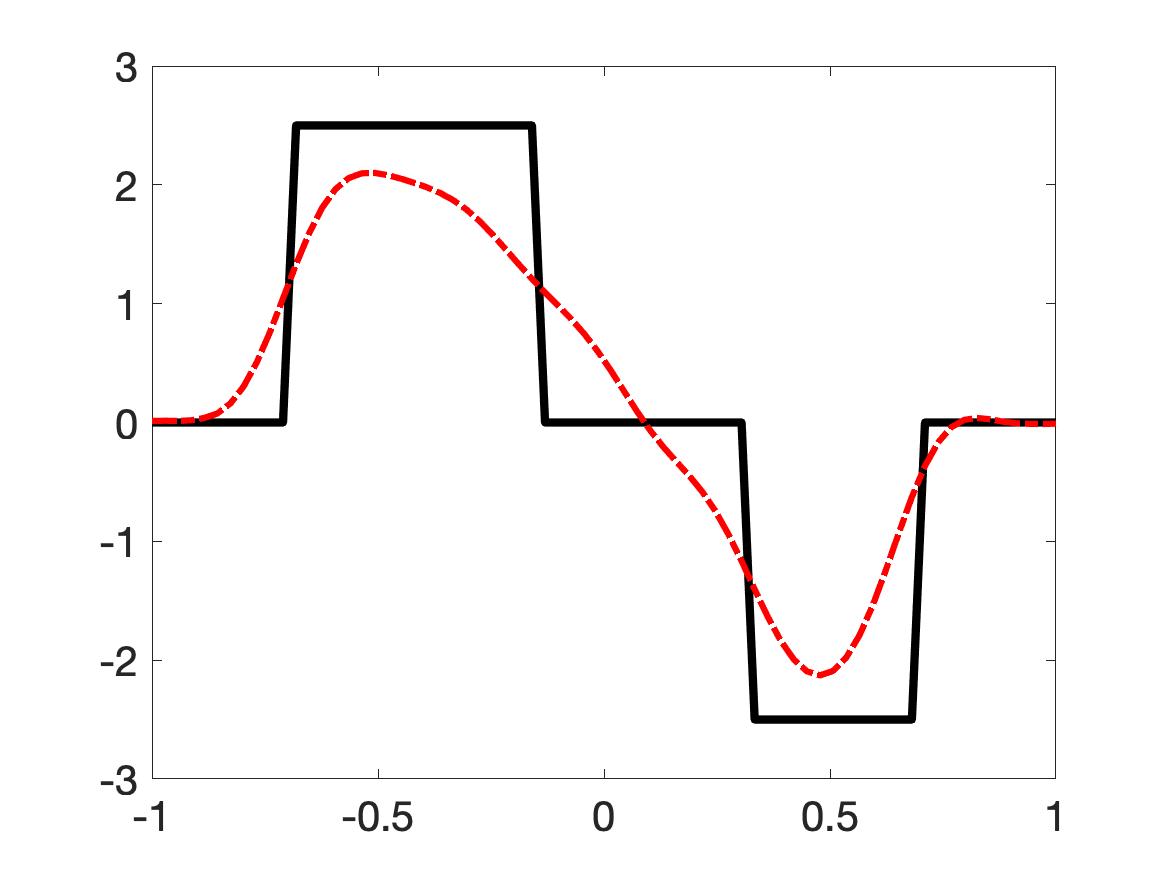}}\hfill 
\subfloat[\label{fig test 4 e}The function $p^{\rm comp}$ computed by Algorithm
\ref{alg} with 80\% noise in the
data]{\includegraphics[width=0.3\textwidth]{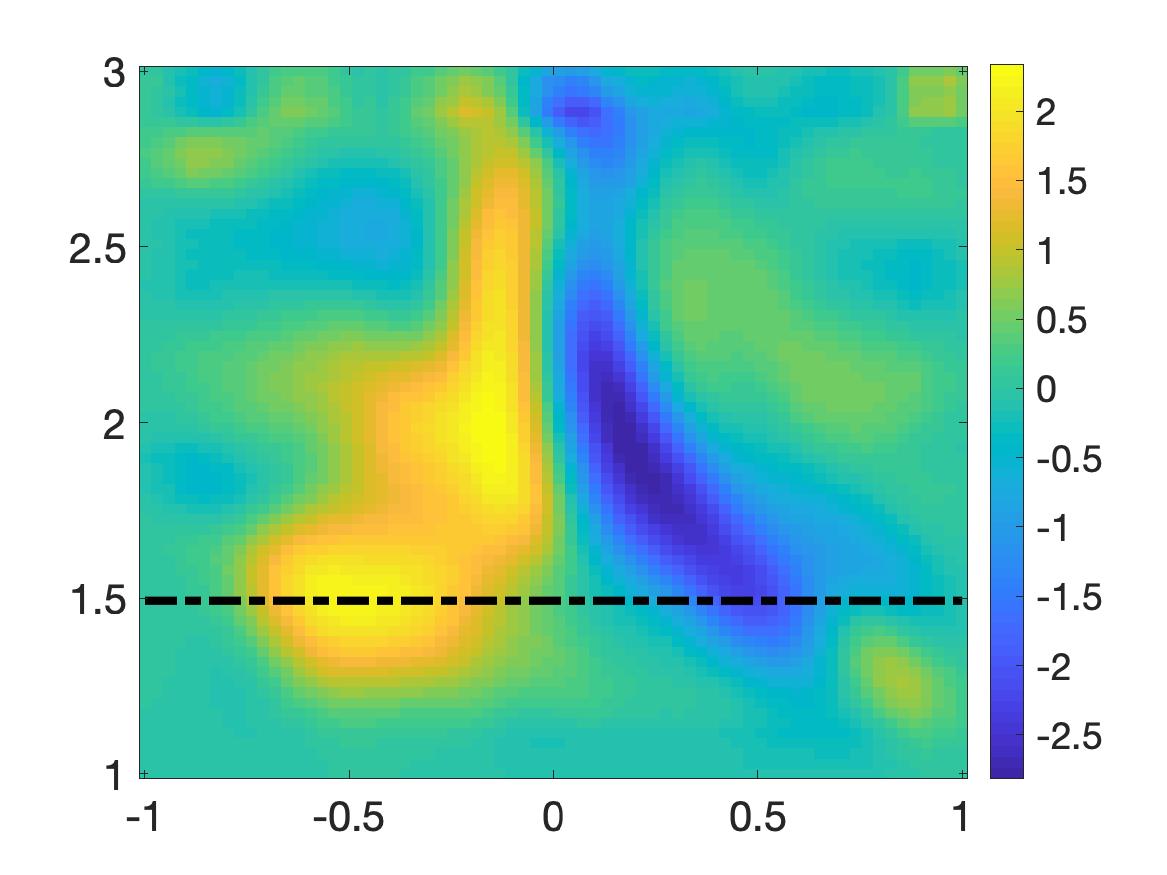}}\hfill 
\subfloat[\label{fig test 4 f}The function $p_{\rm true}$ and $p^{\rm comp}$
with 80\% noise in the data on the set $\{z = 1.5\},$ indicated by a dash-dot line
in (e)]{\includegraphics[width=0.3\textwidth]{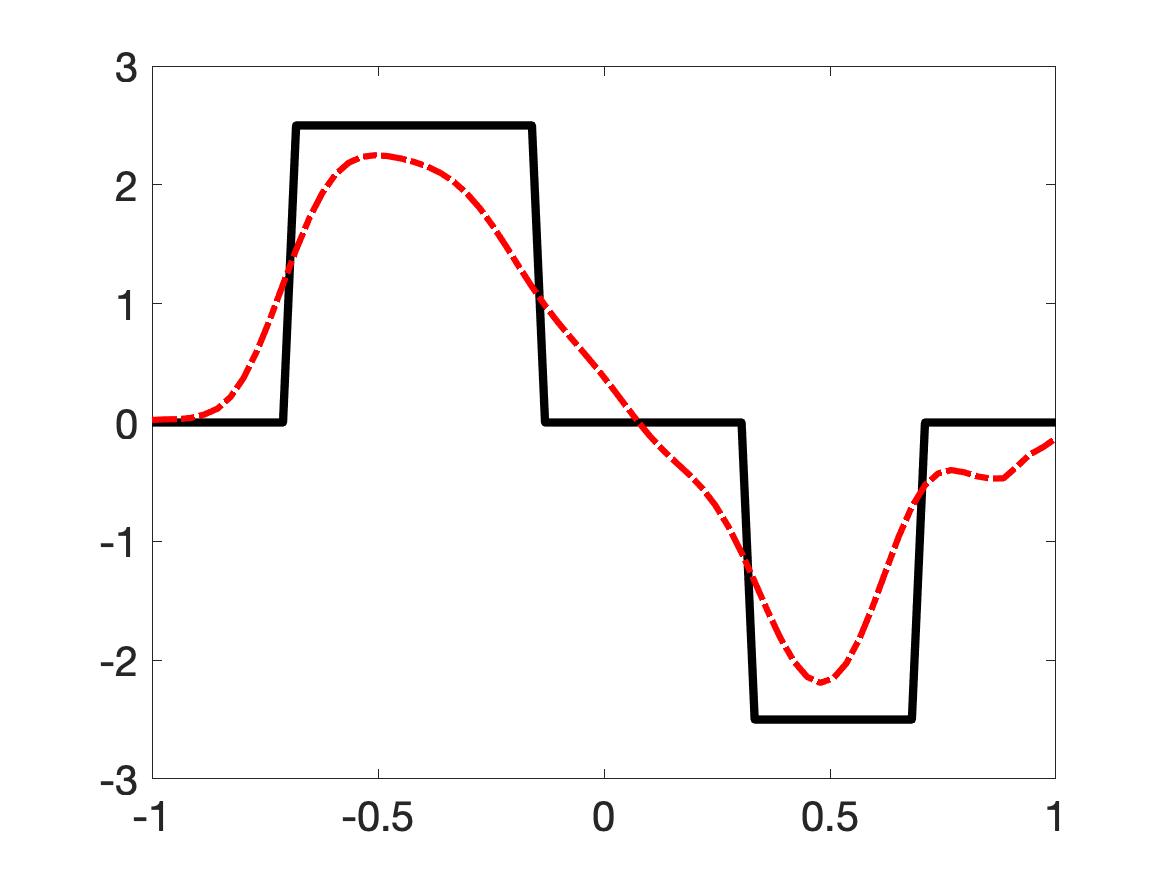}}
\end{center}
\caption{Test 3. The true and reconstructed source functions using Algorithm 
\protect\ref{alg} from noisy data. }
\label{fig test 4}
\end{figure}

\noindent \textbf{Test 3.} We test an interesting and complicated case of
the up-side-down letter $Y$ having both positive and negative values. In
this test, the function $p_{\mathrm{true}}$ is given by 
\begin{equation*}
p^{\mathrm{true}}\left( x,z\right) =\left\{ 
\begin{array}{ll}
2.5 & |x-(z-2)|<0.35,\max \{|x|,|z-2|\}<0.7,z<2, x<0, \\ 
-2.5 & |x+(z-2)|<0.2,\max \{|x|,|z-2|\}<0.7,z<2, x>0, \\ 
2.5 & |x|<0.2,\max \{|x|,|z-2|\}<0.8,z>2, x<0, \\ 
-2.5 & |x|<0.2,\max \{|x|,|z-2|\}<0.8,z>2, x>0.%
\end{array}%
\right.
\end{equation*}%
The background function $c_{0}$ is given by 
\begin{equation*}
c_{0}\left( x,z\right) =\left\{ 
\begin{array}{ll}
1+0.5(x+0.5)^{2}\ln (z) & z>1, \\ 
1 & \mbox{otherwise.}%
\end{array}%
\right.
\end{equation*}%
The numerical results of this test are displayed in Figure \ref{fig test 4}.

It is clear from Figure \ref{fig test 4} that both positive and negative
parts of the function $p\left( x,z\right) $ are successfully identified.
When the noise level $\delta =5\%$, the reconstructed maximal value of the
positive part of $p^{\text{comp}}$ is $2.25$ (relative error 10.0\%) and the
reconstructed minimal value of $p^{\text{comp}}$ of the negative part is $%
-2.74$ (relative error $9.6\%.$) When the noise level is $\delta =80\%$, the
reconstructed maximal value of $p^{\text{comp}}$ of the positive part is $%
2.30$ (relative error 8.0\%) and the reconstructed minimal value of $p^{%
\text{comp}}$ of the negative part is $-2.82$ (relative error $12.8\%.$)

\begin{figure}[h!]
\begin{center}
\subfloat[The function $p^{\rm
true}$]{\includegraphics[width=0.3\textwidth]{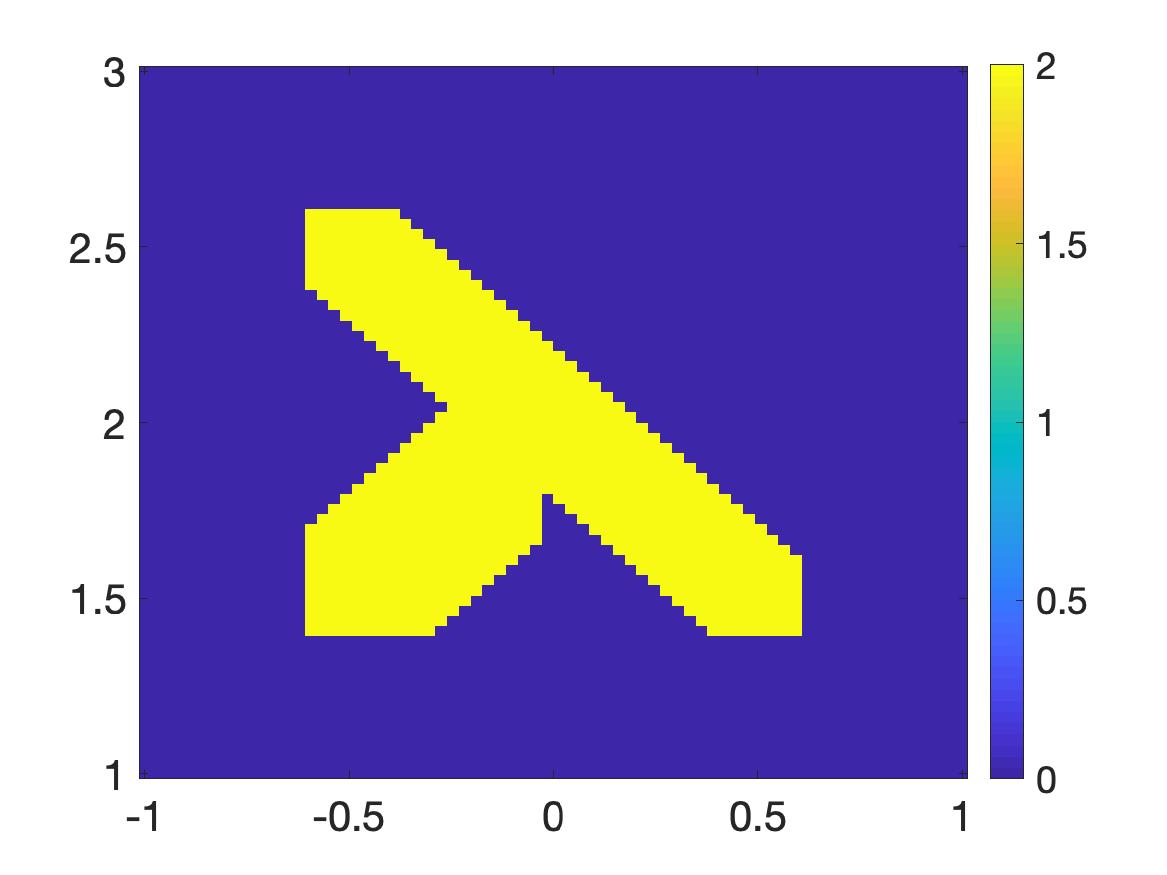}}\hfill 
\subfloat[The function $c_0$ and some geodesic line, generated by the Fast
Marching package in
Matlab]{\includegraphics[width=0.3\textwidth]{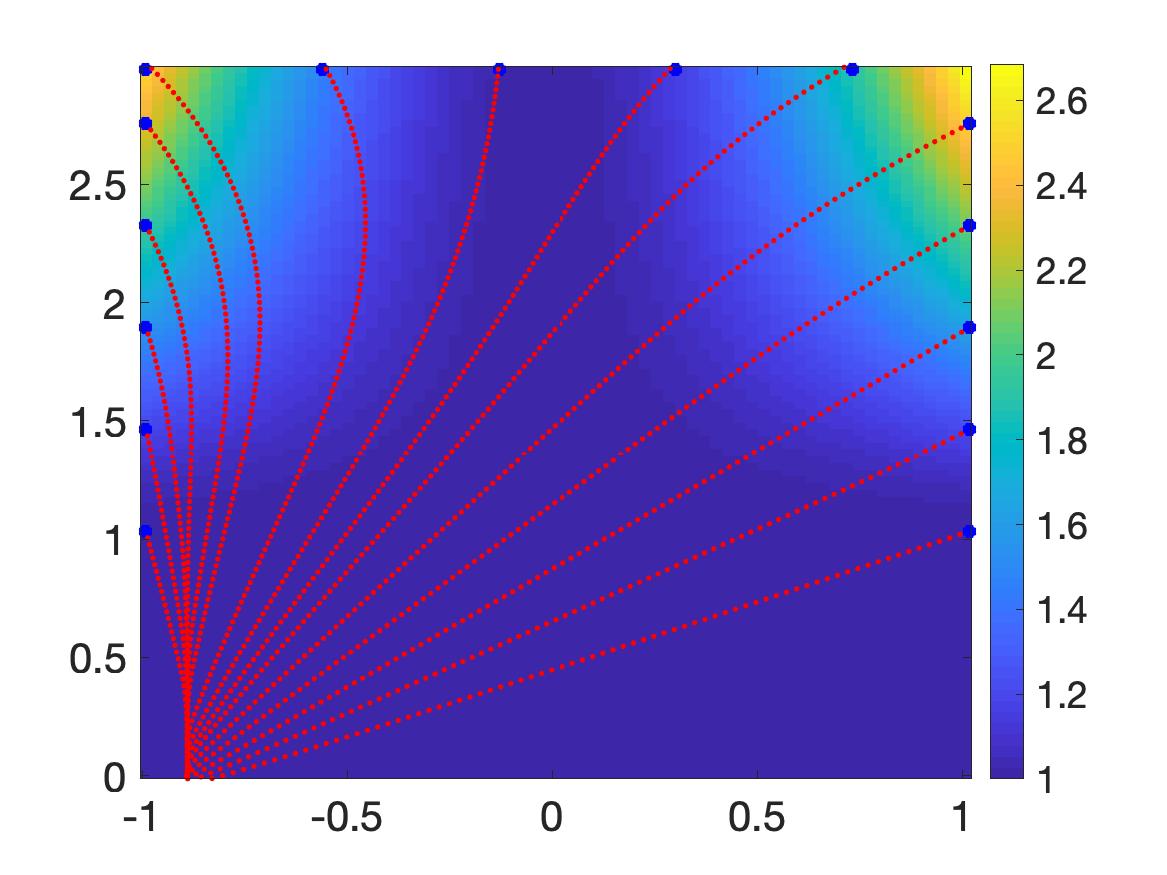}}%
\hfill 
\subfloat[\label{fig test 5 c}The function $p^{\rm comp}$ by
Algorithm \ref{alg} with 5\% noise in the
data]{\includegraphics[width=0.3\textwidth]{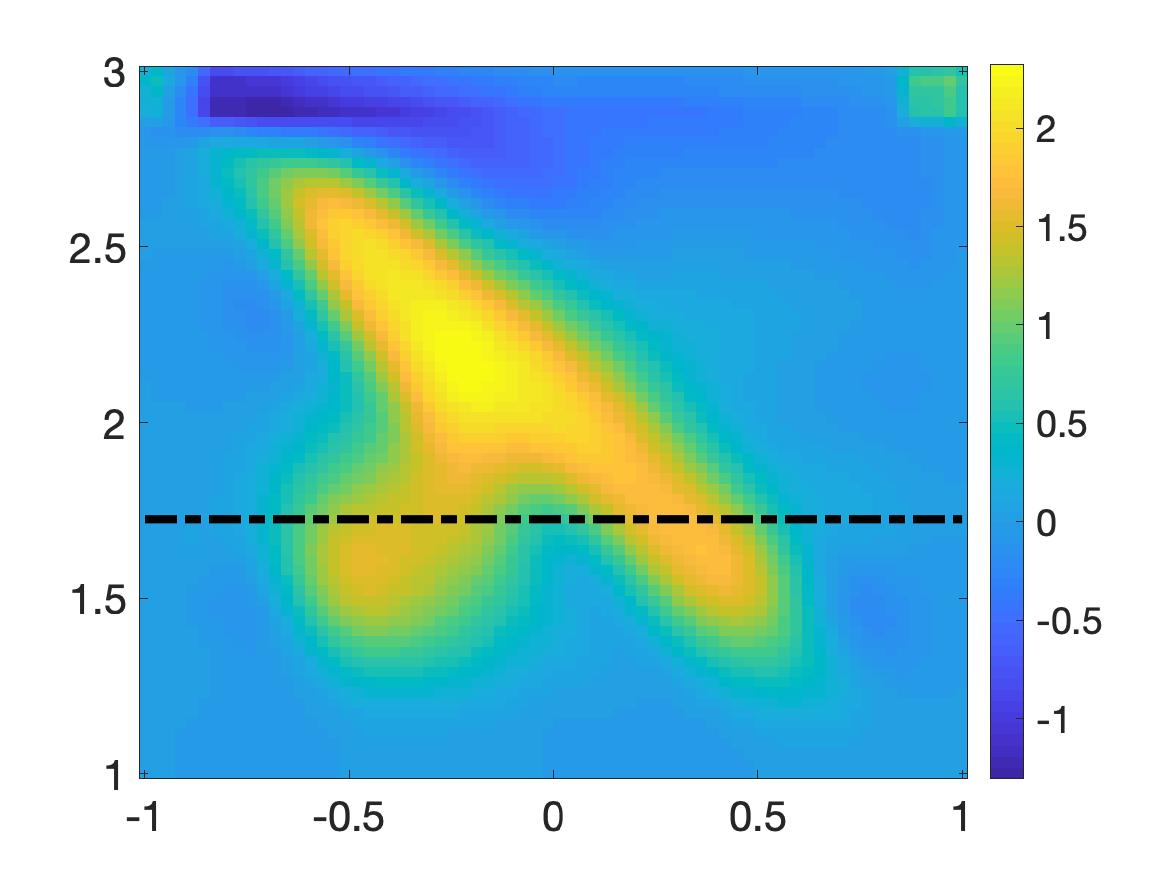}}
\par
\subfloat[\label{fig test 5 d}The function $p^{\rm true}$ and $p^{\rm comp}$
with 5\% noise in the data on the set $\{z = 1.7\},$ indicated by a dash-dot line
in (c)]{\includegraphics[width=0.3\textwidth]{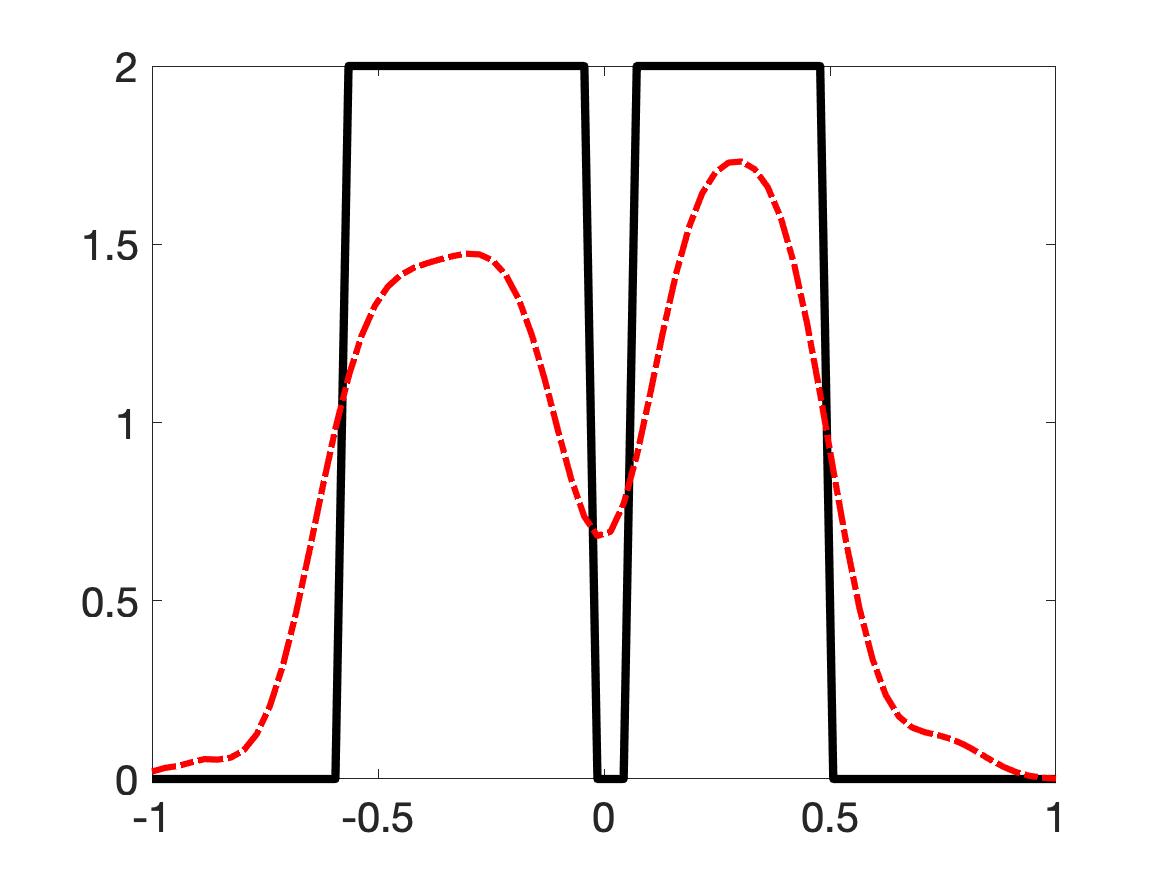}}\hfill 
\subfloat[\label{fig test 5 e}The function $p^{\rm comp}$ computed by Algorithm
\ref{alg} with 80\% noise in the
data]{\includegraphics[width=0.3\textwidth]{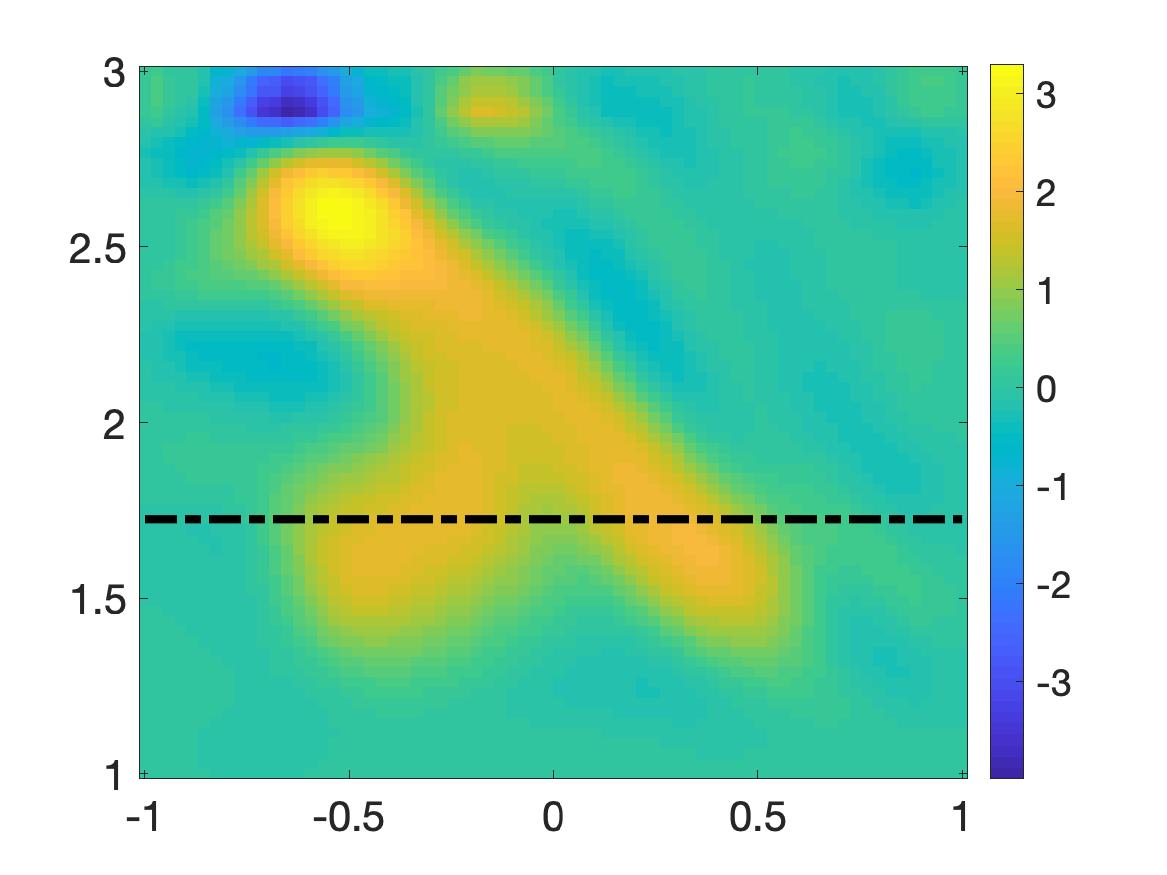}}\hfill 
\subfloat[\label{fig test 5 f}The function $p_{\rm true}$ and $p^{\rm comp}$
with 80\% noise in the data on the set $\{z = 1.7\},$ indicated by a dash-dot line
in (e)]{\includegraphics[width=0.3\textwidth]{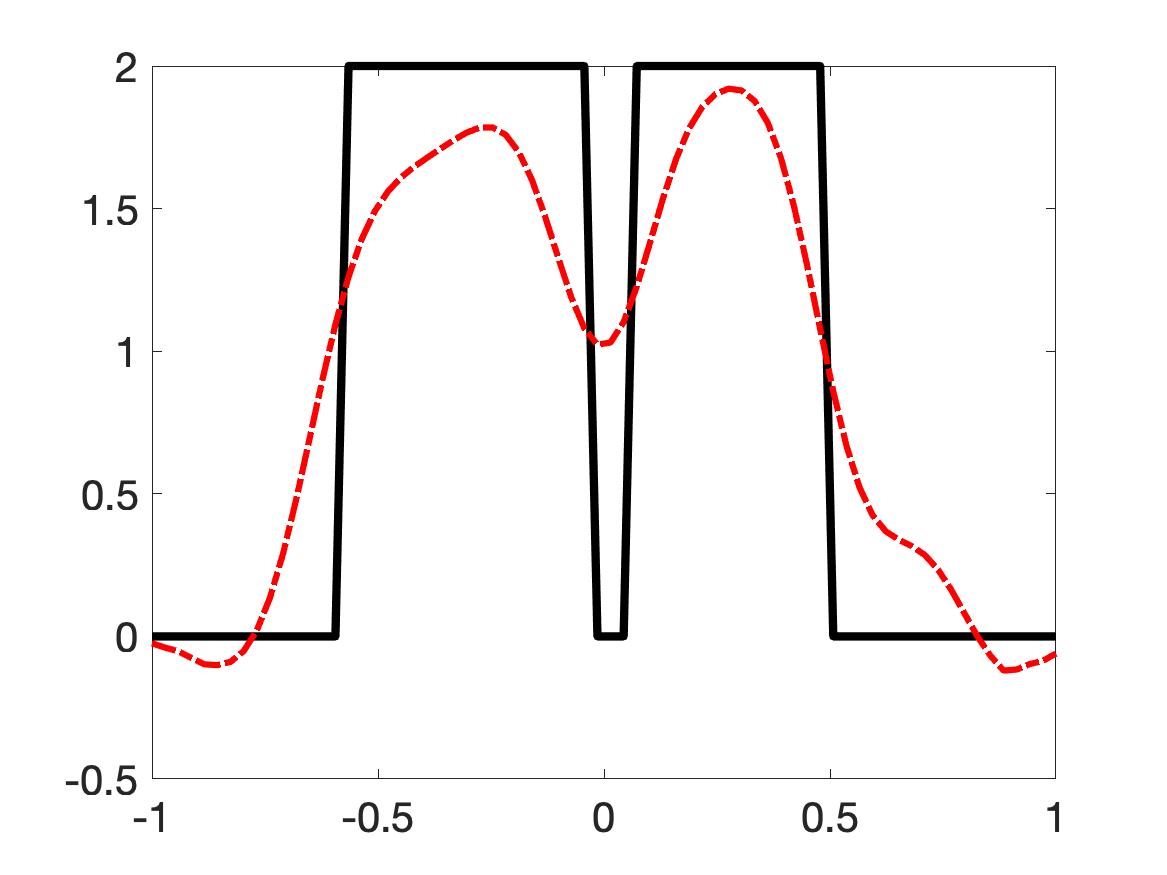}}
\end{center}
\caption{Test 4. The true and reconstructed source functions using Algorithm 
\protect\ref{alg} from noisy data. }
\label{fig test 5}
\end{figure}

\noindent \textbf{Test 4.} In this test, we reconstruct the letter $\lambda $%
. The function $p^{\mathrm{true}}$ is given by 
\begin{equation*}
p^{\mathrm{true}}\left( x,z\right) =\left\{ 
\begin{array}{ll}
2 & |x-(z-2)|<0.325,\max \{|x|,|z-2|\}<0.7\mbox{ and }x<-0.03, \\ 
2 & |x+(z-2)|<0.2\mbox{ and }\max \{|x|,|z-2|\}<0.7, \\ 
0 & \mbox{otherwise.}%
\end{array}%
\right.
\end{equation*}%
In this test, we chose $c_{0}$ as 
\begin{equation*}
c_{0}\left( x,z\right) =\left\{ 
\begin{array}{ll}
1+x^{2}\ln (z) & z>1, \\ 
1 & \mbox{otherwise.}%
\end{array}%
\right.
\end{equation*}%
The numerical results of this test are displayed in Figure \ref{fig test 5}.

The letter $\lambda $ and the values of the function $p^{\mathrm{true}}$ are
successfully reconstructed. The computed position of $\lambda $ is a quite
accurate one, see Figures \ref{fig test 5 d} and \ref{fig test 5 f}. When
the noise level $\delta =5\%,$ the computed maximal value of $p^{\text{comp}%
} $ is 2.31 (relative error 15.5\%). When the noise level $\delta =100\%,$
the computed maximal value of $p^{\text{comp}}$ is 3.27 (relative error
63.5\%).

\section{Concluding Remarks}

\label{sec conclusion} In this paper, we have developed a convergent
numerical method of the solution of the linearized Travel Time Tomography
Problem with non-redundant incomplete data. A good accuracy of numerical
results with 5\% noise in the data is demonstrated for rather complicated
functions to be imaged. It is quite surprising that an acceptable accuracy
of computational results is observed even for very high level of noise in
the data varying between 30\% and 120\%.

%\begin{center}
%\textbf{Acknowledgment}
%\end{center}
%
%This work was supported by US Army Research Laboratory and US Army Research
%Office grant W911NF-19-1-0044.

%\bibliographystyle{siamplain}
%\bibliography{../../../../../mybib}

\end{document}